\documentclass[11pt]{article}
\usepackage{amscd, amssymb, latexsym, amsmath}

\newtheorem{thm}{Theorem}

\newtheorem{pro}{Proposition}
\newtheorem{lem}{Lemma}
\newtheorem{dfn}{Definition}

\newtheorem{rem}{Remark}

\newenvironment{proof}
{\noindent {\em Proof.}} {\hfill q.e.d.}

\numberwithin{thm}{section} \numberwithin{cor}{section}
\numberwithin{pro}{section} \numberwithin{lem}{section}
\numberwithin{dfn}{section}
\numberwithin{rem}{section}
\numberwithin{equation}{section}

\newcommand{\R}{\mathbb R}

\newcommand{\barbla}{\overline{\nabla}}
\newcommand{\baromega}{\overline{\omega}}
\newcommand{\pk}{\partial_k}
\newcommand{\pl}{\partial_l}
\newcommand{\pai}{\partial_i}
\newcommand{\pj}{\partial_j}
\begin{document}
\title{Mean Curvature Flow of Surfaces
\\ in Einstein Four-Manifolds}
\author{Mu-Tao Wang}
\date{September 21, 2000, revised November 28, 2000}
\maketitle

\begin{abstract}
Let $\Sigma$ be a compact oriented surface immersed in a four
dimensional K\"ahler-Einstein manifold $(M, \omega)$. We consider
the evolution of $\Sigma$ in the direction of its mean curvature
vector. It is proved that being symplectic is preserved along the
flow and the flow does not develop type I singularity.
 When $M$ has two parallel K\"ahler forms
$\omega'$ and $\omega''$ that determine different
orientations and
$\Sigma$ is symplectic with respect to
both $\omega'$ and $\omega''$, we prove the mean curvature flow of
$\Sigma$ exists smoothly for all time.
In the positive curvature case, the
flow indeed converges at infinity.
\end{abstract}

\section{Introduction}
Let $(M,g)$ be a Riemannian manifold
and let $\alpha$ be a calibrating $k$-form on $M$
i.e. $d\alpha=0$ and comass$(\alpha)=1$. In
this article, we shall assume additionally
$\alpha$ is parallel. This in particular
implies $M$ is of special holonomy.

A $k$-dimensional submanifold is said to be
calibrated by $\alpha$ if the restriction
of $\alpha$ gives the volume form of the
submanifold. A simple
application of Stoke's theorem shows a
calibrated submanifold minimizes the volume
functional in its homology class. To produce
 a calibrated submanifold, it is
 thus natural to consider the
gradient flow of the volume functional. By the
first variation formula of volume, this is
equivalent to evolving a submanifold $\Sigma_0$
in the
direction of its mean curvature vector. To make it
precise, the mean curvature flow is the
solution of the following system of parabolic
equations.

\[\frac{dF}{dt}(x,t)=H(x,t)
\]

\noindent where $F:\Sigma\times [0,T) \rightarrow M$ is
a one parameter family of immersions $F_t(\cdot)
=F(\cdot, t)$ of
 $\Sigma$ into $M$.
$H(x,t)$ is the mean curvature vector of
$F_t(\Sigma)$ at $F_t(x)$. We say $F$ is the mean curvature
flow of the immersed submanifold $F_0(\Sigma)$.
For a fixed $t$, the submanifold $F_t(\Sigma)$
is denoted by $\Sigma_t$.
If we assume $M=\R^n$. In terms of coordinate
$x^1,\cdots ,x^k$ on $\Sigma$, the mean
curvature flow is the following system of
parabolic equations
\[F=F^A(x^1,\cdots, x^k, t),\,\, A=1,\cdots n\]
\[\frac{\partial F^A}{\partial t}
=\sum_{i, j, B}g^{ij}\,P^A_B\,\frac{\partial^2 F^B}{\partial
x^i\partial x^j}, \,\, A=1,\cdots n \] where $g^{ij}$ is the
inverse matrix to $g_{ij}=\frac{\partial F^A}{\partial x^i}
\frac{\partial F^A}{\partial x^j}$ and $P^A_B
=\delta^A_B-g^{kl}\frac{\partial F^A} {\partial x^k}
\frac{\partial F^B}{\partial x^l}$ is the projection to the normal
part.

The mean curvature flow of hypersurfaces has been studied
extensively in the last decade. In this case, the mean curvature
$H$ is essentially a scalar function and the positivity of $H$ is
preserved along the flow. Very little is known in higher
codimension except for the curve flows.

This article considers the next simplest higher codimension mean
curvature flow, namely a surface flows in a four dimensional
manifold. We impose a positivity condition on the initial
submanifold. An oriented submanifold $\Sigma$ is said to be almost
calibrated by $\alpha$ if $*\alpha>0$ where $*$ is the Hodge star
operator on $\Sigma$.

The following question arises naturally. Can an almost calibrated
submanifold be deformed to a calibrated one along the mean
curvature flow? We study this question in the case when $M$ is a
four-dimensional Einstein manifold and $\Sigma_0$ is almost
calibrated by a parallel calibrating form. When $M$ is a
K\"ahler-Einstein surface and the calibrating form is the K\"ahler
form, an almost calibrated surface is a symplectic curve with the
induced symplectic structure. A calibrated submanifold in this
case is a holomorphic curve.

We use blow up analysis to characterize the singularities of
mean curvature flow of symplectic surfaces.
It turns out they are all so-called type II
singularities.

\vskip 10pt
\noindent {\bf Theorem A}
\hskip 3pt{\it Let $M$ be a four-dimensional
K\"ahler-Einstein manifold, then a symplectic
surface remains symplectic along the
mean curvature flow and the flow does not
develop any type I singularities.}
\vskip 10pt

When $M$ is
locally a product and the initial surface
is almost calibrated by two calibrating forms, we prove the following
long time existence theorem.

\vskip 10pt \noindent {\bf Theorem B} \hskip 3pt{\it Let $M$ be an
oriented four-dimensional Einstein manifold with two parallel
calibrating forms $\omega', \omega''$ such that $\omega'$ is
self-dual and $\omega''$ is anti-self-dual. If $\Sigma$ is a
compact oriented surface immersed in $M$ such that $*\omega',
*\omega''>0$ on $\Sigma$. Then the mean curvature flow of $\Sigma
$ exists smoothly for all time.} \vskip 10pt

We remark that the assumption implies $M$ is locally a product of
two surfaces. As for convergence at infinity, we prove the
following theorem in the non-negative curvature case.

\vskip 10pt
\noindent {\bf Theorem C}
\hskip 3pt{\it Under the same assumption as
in Theorem B. When $M$ has non-negative
curvature, there exists a constant
$1>\epsilon>0$ such that if $\Sigma$ is a compact oriented surface
immersed in $M$ with $*\omega',
*\omega''>1-\epsilon $ on $\Sigma$, the mean
curvature flow of $\Sigma $ converges
smoothly to a totally geodesic
surface at infinity.}
\vskip 10pt

This is proved by an uniform estimate
of the norm of the second fundamental form.

When $M=S^2\times S^2$, the combination of
Theorem B and C yields

\vskip 10pt \noindent {\bf Theorem D} \hskip 3pt{\it Let $M=(S^2,
\omega_1)\times (S^2,\omega_2)$. If $\Sigma$ is a compact oriented
surface embedded in $M$ such that $*\omega_1>|*\omega_2|$. Then
the mean curvature flow of $\Sigma$ exists for all time and
converges smoothly to an $S^2\times \{p\}$.} \vskip 10pt

This theorem in  particular applies to the graph of
maps between two Riemann surfaces. Namely,
let $f:(\Sigma_1, \omega_1)\rightarrow
(\Sigma_2, \omega_2)$ be a map between Riemann
surfaces of the same constant curvature and
$\omega_i$ is the volume form of $\Sigma_i$.
We consider the product $M=\Sigma_1\times
\Sigma_2$ and let $\omega'=\omega_1+\omega_2$
and $\omega''=\omega_1-\omega_2$. If the
Jacobian of $f$ is less than one, then
we have $*\omega'>0$ and $*\omega''>0$ on
the graph of $f$. Therefore this formulation
gives a natural way to deform the map
$f$ to a constant map.

\vskip 5pt
\noindent {\bf Corollary D}
\hskip 3pt{\it Any smooth map between two-spheres
with Jacobian less than one deforms to
a constant map through the mean curvature
flow of the graph.}
\vskip 8pt

The article is organized as the followings.
In \S 2, the parabolic equation satisfied
by a general parallel form along the mean
curvature flow is derived. \S 3
discusses general calibrating two-forms in
a four dimensional space. \S 4 computes
the equation satisfied by a K\"ahler
form along the mean curvature flow. \S 5 studies the singularities of
 mean curvature flow of symplectic surfaces
 and proves Theorem A.
 \S 6 concerns
long time existence and Theorem B is proved
there.  Convergence
at infinity is discussed in \S 7. Theorem C
is proved at the end of this section. \S 8
discusses applications in the positive curvature
case and proves Theorem D.

This project starts in the fall of 1998 in an attempt to answer
Professor S.-T. Yau's question " how to deform a symplectic
submanifold to a holomorphic one". Theorem A, in particular the
result "symplectic remains symplectic" and the exclusion of type I
singularity, was obtained in the summer of 1999. It has been
presented in the geometry seminars at Stanford, U. C. Berkeley, U.
C. Santa Cruz and U of Minnesota between February 2000 and May
2000. I would like to thank Professor R. Schoen and Professor
S.-T. Yau for their constant encouragement and invaluable advice.
I also have benefitted greatly from the many discussion that I
have with Professor G. Huisken, Professor L. Simon and Professor
B. White.

\section{Evolution equations of parallel
forms}

Let $F:\Sigma^2 \rightarrow M^{4}$ be an
isometric immersion of an orientable surface
into a four-dimensional Riemannian
manifold. We fixed an orientation on $\Sigma$.
The
restriction of the tangent bundle of $M$ to $\Sigma$ splits
as the direct sum of the tangent bundle of $\Sigma$
and the normal bundle.
\[
TM|_\Sigma=T\Sigma\oplus N\Sigma
\]

The Levi-Civita connection on $M$ induces
a connection on $T\Sigma$. We denote
the connection on $M$ by $\barbla $ and
 the induced connection on $T\Sigma$ by $\nabla$.
Therefore,

\[
\nabla_X Y=(\barbla_X Y)^T
\]
for any tangent vector fields $X, Y$.
Here $(\cdot)^T$ denotes the projection from
$TM$ onto $T\Sigma$ and $(\cdot)^N$ shall denote the
projection onto $N\Sigma$.

The second fundamental form $A:T\Sigma
\times T\Sigma \mapsto N\Sigma$ is
defined by $A(X,Y)= (\barbla_{X} Y)^N$. We also
define $B:T\Sigma\times N\Sigma\mapsto T\Sigma$
by $B(X, N)=(\barbla_X N)^T$. The relation
between $A$ and $B$ is
\[<A(X,Y), N>=-<Y, B(X,N)>\]

Notice that
 we have identified $X\in T\Sigma$ with
 $F_*(X)\in TM$.

Fix a point $p\in \Sigma$. Let $\{x^i\}$ be
a normal coordinate system
for $\Sigma$ at $p$ and $\{y^A\}$ a
normal coordinate system for $M$ at $F(p)$. We denote
$\frac{\partial}{\partial x^i}$ by $\partial_i$
and identify it with $\frac{\partial F}{\partial x_i}$.
The induced metric on $\Sigma$ is given
by $g_{kl}=<\pk,\pl>$. The mean curvature
vector along $\Sigma$ is the trace of $A$, i.e $H=g^{kl}
A(\pk,\pl)$.

Let $\overline{\omega}$ be a parallel two
form on $M$ and $\omega=F^*\overline{\omega}$
be the pull-back of $\omega$ on $\Sigma$. We
first compute the rough Laplacian of $\omega$
on $\Sigma$.

\[\Delta \omega= g^{kl}\nabla_{\pk}\nabla_{\pl}
\omega
\]

\begin{lem}
\begin{equation}
\begin{split}
(\Delta\omega)(X,Y)
&=\baromega((\barbla_X H)^N,Y)
-\baromega((\barbla_Y H)^N,X)\\
&-g^{kl}\baromega((K(\pk, X)\pl)^N,Y)
+g^{kl}\baromega((K(\pk, Y)\pl)^N,X)\\
&+g^{kl}\baromega(B(\pk, A(\pl, X)), Y)
-g^{kl}\baromega(B(\pk, A(\pl, Y)), X)\\
&+2g^{kl}\baromega(A(\pk,X), A(\pl,Y))
\end{split}
\end{equation}
where $K(X,Y)Z=-\barbla_X \barbla_Y Z
+\barbla_Y \barbla_X Z-\barbla_{[X,Y]}Z$
is the curvature operator of $M$. Notice
that $<K(X,Y)X,Y>>0$ if $M$ has positive sectional
curvature.
\end{lem}

\begin{proof}
Since both sides are tensors, we calculate at the point $p$ using
a normal coordinate system. Therefore $g_{kl}=\delta_{kl}$ and all
connection terms vanish. Now
\[(\Delta\omega)(\pai,\pj)
=\pk[\pk(\omega(\pai,\pj))-\omega(\nabla_{\pk}
\pai, \pj)-\omega(\pai, \nabla_{\pk}\pj)]\]

The term in the bracket is
\[
\begin{split}
&\pk(\omega(\pai,\pj))-\omega(\nabla_{\pk}
\pai, \pj)-\omega(\pai, \nabla_{\pk}\pj)\\
&=(\barbla_{\pk}\baromega)(\pai,\pj)
+\baromega(\barbla_{\pk} \pai, \pj)+\baromega
(\pai, \barbla_{\pk} \pj)-
\baromega(\nabla_{\pk}\pai,\pj)
-\baromega(\pai, \nabla_{\pk} \pj)\\
&=\baromega(A(\pk,\pai), \pj)
+\baromega(\pai, A(\pk, \pj))
\end{split}
\]
where we have used the fact that $\baromega$ is
parallel and $\barbla_{\pk} \pai-\nabla_{\pk}
\pai=A(\pk,\pai)$.

Therefore
\begin{equation}\label{lapom}
\Delta\omega(\pai,\pj)
=\pk[\baromega(A(\pk,\pai), \pj)
+\baromega(\pai, A(\pk, \pj))]
\end{equation}

Use Leibnitz rule and the parallelity of
$\baromega$ again.
\[
\begin{split}
&\pk(\baromega(A(\pk,\pai), \pj)\\
&=\baromega(\barbla_{\pk} A(\pk, \pai), \pj)
+\baromega(A(\pk, \pai), \barbla_{\pk}\pj)\\
&=\baromega((\barbla_{\pk} A(\pk, \pai))^T
+(\barbla_{\pk} A(\pk, \pai))^N, \pj)
+\baromega(A(\pk, \pai), \nabla_{\pk}\pj
+A(\pk,\pj))\\
&=\baromega(B({\pk}, A(\pk, \pai)),\pj)
+\baromega ( A(\pk, \pai), A(\pk, \pj))
+\baromega((\barbla_{\pk} A(\pk, \pai))^N, \pj)
\end{split}
\]
where we have used $\nabla_{\pk} \pj =0$
at the point $p$ in normal coordinates.

\[
\begin{split}
&\baromega((\barbla_{\pk} A(\pk, \pai))^N, \pj)\\
&=\baromega((\barbla_{\pk}  \barbla_{\pai}\pk)^N, \pj)
-\baromega((\barbla_{\pk} \nabla_{\pai} \pk)^N, \pj)\\
&=\baromega((-K(\pk, \pai)\pk,
+\barbla_{\pai} \barbla_{\pk} \pk)^N, \pj )
-\baromega(B(\pk, \nabla_{\pai}\pk),\pj)\\
&=-\baromega((K(\pk, \pai)\pk)^N, \pj)
+\baromega((\barbla_{\pai}H)^N+
 (\barbla_{\pai} \nabla_{\pk}\pk)^N, \pj )\\
&=-\baromega((K(\pk, \pai)\pk)^N, \pj)
+\baromega((\barbla_{\pai}H)^N, \pj)+
 \baromega(B({\pai}, \nabla_{\pk}\pk)), \pj )\\
\end{split}
\]
The last term vanishes in normal coordinates.

Thus we have proved

\[
\begin{split}
\pk(\baromega(A(\pk,\pai), \pj)
&=\baromega(B(\pk, A(\pk,
\pai)),\pj)
+\baromega(A(\pk, \pai), A(\pk, \pj))\\
&-\baromega((K(\pk, \pai)\pk)^N, \pj)
+\baromega((\barbla_{\pai}H)^N, \pj)
\end{split}
\]

Plug this equation back into equation (\ref{lapom}) and
anti-symmetrize $i$, $j$ and the lemma is proved.
\end{proof}

Let's represent the fixed orientation on $\Sigma$ by
a two-form $d\mu$. Let $F:\Sigma\times [0, T) \rightarrow M$ be
the mean curvature flow of $\Sigma$.
The immersion $F_t$ induces a pull-back metric
$g_t$ on $\Sigma$. The volume form of $g_t$ is
denoted by $d\mu_t=\sqrt{\det g_t}\,d\mu$.

Now we consider the evolution equation of
$\omega_t= F_t^*(\baromega) $. This is
a family of time-dependent two forms on the
fixed surface $\Sigma$. Let the one-form
$\alpha_t$ be defined by  $\alpha_t(X)=\baromega
(H_t, X)$.

\begin{lem}
Along the mean curvature flow
\[\frac{d}{dt}\omega_t= d\alpha_t\]
For any vector field $X, Y\in T\Sigma$,
\[\frac{d}{dt}\omega_t(X,Y)
=\baromega((\barbla_{X} H)^N , Y)
+\baromega(X, (\barbla_{Y}H)^N)
+\baromega(B(X, H), Y)
+\baromega(X, B(Y, H))
\]
\end{lem}

\begin{proof}

\[\frac{d}{dt}\omega_t(\pai,\pj)
=\baromega(\barbla_{H} \pai, \pj)
+\baromega(\pai, \barbla_{H}\pj)
 =\baromega(\barbla_{\pai} H , \pj)
+\baromega(\pai, \barbla_{\pj}H)
\]

By definition $\barbla_{\pai}H=(\barbla_{\pai}H)^N
+B(\pai, H)$.

On the other hand

\[d\alpha_t (\pai, \pj)
=\pai(\baromega(H, \pj))-\pj(\baromega(H, \pai))
=\baromega(\barbla_{\pai} H, \pj)
-\baromega(\barbla_{\pj} H, \pai)
\]

\end{proof}

The volume form $d\mu_t$ determines a Hodge
operator $*_t$. Therefore
$*_t \omega_t$ becomes a time-dependent function on $\Sigma$.

\begin{pro}\label{star}
Let $\baromega$ be a parallel two-form on $M$.
$F_t:\Sigma\mapsto M$ be the $t$ slice
of a  mean curvature flow
and $\omega_t=F_t^*(\baromega)$ be the
pull-back form on $\Sigma$. Then $\eta_t
=*_t \omega_t$ satisfies the following
parabolic equation.

\[
\begin{split}
\frac{d}{dt}\eta_t&=(\Delta_t \eta_t)+
|A|^2\,\eta_t\\
&-2\baromega(A(e_k,e_1), A(e_k,e_2))
+\baromega((K(e_k, e_1)e_k)^N,e_2)
-\baromega((K(e_k, e_2)e_k)^N,e_1)\\
\end{split}
\]
where $|A|$ is the norm of the second
fundamental form, $|A|^2
=g^{ij}g^{kl}<A(\pai,\pk), A(\pj, \pl)>$ and $\{e_1, e_2\}$ any
orthonormal basis with respect to $g_t$.
\end{pro}

\begin{proof}
Combine the previous two lemma, we get

\begin{equation}\label{form}
\begin{split}
\frac{d}{dt}\omega_t (X,Y)&=(\Delta_t \omega_t)
(X,Y)+\baromega(B(X,H), Y)
+\baromega(X, B(Y,H))\\
&+g^{kl}\baromega((K(\pk, X)\pl)^N,Y)
-g^{kl}\baromega((K(\pk, Y)\pl)^N,X)\\
&-g^{kl}\baromega(B(\pk, A(\pl, X)), Y)
+g^{kl}\baromega(B(\pk, A(\pl, Y)), X)\\
&-2g^{kl}\baromega(A(\pk,X), A(\pl,Y))
\end{split}
\end{equation}

Now $*_t\omega_t=\frac{\omega(\partial_1
,\partial_2)}{\sqrt{\det g_t}}$ where
$\{\partial_1,\partial_2\}$ is a fixed
coordinate system on $\Sigma$ and
$\det g_t$ is the determinant of
$(g_t)_{ij}=<(F_t)_*\pai, (F_t)_*\pj>$.

It is easy to compute
\[\frac{d}{dt}{\sqrt{\det g_t}}=-|H|^2 {\sqrt{\det g_t}}
\]
where $|H|$ is the norm of the mean curvature
vector.

Thus

\[\frac{d}{dt}*_t\omega_t
=\frac{1}{\sqrt{\det g_t}}\frac{d}{dt}
\omega_t(\partial_1,\partial_2)
+|H|^2 *_t\omega_t
\]

Now we use equation (\ref{form})
with $X=\partial_1$ and $Y=\partial_2$.

The first term is $\frac{1}{\sqrt{\det g_t}}
(\Delta_t\omega_t)(\partial_1,\partial_2)=
*_t \Delta_t \omega_t=\Delta_t *_t\omega_t$
because the Hodge $*_t$ operator is parallel.

For other terms we can take any orthonormal
basis $\{e_1, e_2\}$ with respect to the metric
$g_t$ to calculate.

It is not hard to see

\[
\begin{split}
&\baromega(B(e_1,H), e_2)
-\baromega(B(e_2,H),e_1)\\
&=*_t\omega_t
(< B(e_1, H),e_1>+< B(e_2, H),e_2>)\\
&=-*_t\omega_t(<A(e_1, e_1), H>+A(e_2,e_2), H>\\
&=-*_t\omega_t |H|^2
\end{split}
\]

Likewise,

\[
\begin{split}
&\baromega(B(e_k, A(e_k, e_1)), e_2)
-\baromega(B(e_k, A(e_k, e_2)), e_1)\\
&=*_t\omega_t (<B(e_k, A(e_k, e_1)), e_1>+
<B(e_k, A(e_k, e_2)), e_2>)\\
&=-*_t\omega_t(<A(e_k, e_1),A(e_k, e_1)>
+<A(e_k, e_2),A(e_k, e_2)>)
\end{split}
 \]

\end{proof}
\section{Calibrating two-forms
in four-dimensional spaces}

Let $V \backsimeq \R^4$ be an inner product space
and $\alpha\in \wedge^2 V^*$ a two form. We shall
 use the inner product to identify $V$ and $V^*$ and
this induces inner product on all $\wedge^k V^*$. First let's
recall the definition of comass of $\alpha$,
\[\text{comass}(\alpha) =\max_{x\in G(2,V)}
\alpha (x)
\]
where $G(2,V)$ is the Grassmanian of all
two-planes in $V$. $G(2,V)$ can be described
by

\[G(2,V)=\{x\in \wedge^2 V, x\wedge x=0\,
\text{and}\, |x|^2=1\}\]

Now we fix an orientation $\nu\in \wedge^4 V^*$
and normalize so that $|\nu|=1$. Given any
orthonormal basis $\{e_1, e_2, e_3, e_4\}$
for $V$ such that $\nu(e_1, e_2, e_3, e_4)=1$,
the following two-forms give an orthonormal
basis for $\wedge^2 V^*$.

\[
\begin{matrix}
&\alpha_1={1\over \sqrt{2}}(e_1^*\wedge
e_2^*+e^*_3\wedge e_4^*)
&\beta_1={1\over \sqrt{2}}(e_1^*\wedge
e_2^*-e^*_3\wedge e_4^*)\\
&\alpha_2={1\over \sqrt{2}}(e_1^*\wedge
e_3^*-e^*_2\wedge e_4^*)
&\beta_2={1\over \sqrt{2}}(e_1^*\wedge
e_3^*+e^*_2\wedge e_4^*)\\
&\alpha_3={1\over \sqrt{2}}(e_1^*\wedge
e_4^*+e^*_2\wedge e_3^*)
&\beta_3={1\over \sqrt{2}}(e_1^*\wedge
e_4^*-e^*_2\wedge e_3^*)\\
\end{matrix}
\]

These forms serve as coordinate functions
on $G(2,V)$, under the identification
\[x\rightarrow (\alpha_i(x), \beta_i(x))\]
An element $x$ in $G(2,V)$ satisfies
$\sum_i (\alpha_i(x))^2=
\sum_i (\beta_i (x))^2={1\over 2}$.
Therefore $G(2,V)\backsimeq S^2(\frac{1}
{\sqrt{2}})\times S^2(\frac{1}
{\sqrt{2}})$.

Now for any given $\alpha \in \wedge^2 V^*$. We identify
$\alpha $ with an element $K$ in $End(V)$
by $\alpha(X, Y)=<K(X), Y>$.
Since $\sqrt{-1}
\,K$ is Hermitian symmetric and purely
imaginary, it
has real eigenvalues $\pm \lambda_1, \pm \lambda_2$.
We can require
$\alpha \wedge \alpha =\lambda_1 \lambda_2
\nu$. $\lambda_1 \lambda_2$ is
actually the
Pfaffian of $\alpha$ and $\det \, K
=(\lambda_1 \lambda_2)^2$. A form is
self-dual (anti-self-dual ) if
$\lambda_1 \lambda_2=1(-1)$.

\begin{lem}
\[\text{comass} (\alpha)
=\max \{|\lambda_1|, |\lambda_2|\}\]
\end{lem}

\begin{proof}
We can find an orthonormal basis $\{e_1,e_2,
e_3, e_4\}$ with $\nu(e_1, e_2, e_3, e_4)=1$
such that $\alpha =\lambda_1 e_1^*\wedge
e_2^*+\lambda_2 e_3^*\wedge e_4^*$. In terms
of the self-dual and anti-self-dual bases
associated with $\{e_1,e_2,
e_3, e_4\}$.

\[\alpha ={1\over\sqrt{2}}(\lambda_1+\lambda_2)
\alpha_1+{1\over\sqrt{2}}(\lambda_1-\lambda_2)
\beta_1
\]

Therefore
\[
\begin{split}
\alpha(x)&= {1\over\sqrt{2}}(\lambda_1+\lambda_2)
\alpha_1(x)+{1\over\sqrt{2}}(\lambda_1-\lambda_2)
\beta_1(x)\\
& \leq
{1\over 2}(|\lambda_1+\lambda_2|
+|\lambda_1-\lambda_2|)\\
&=\max\{\lambda_1, \lambda_2\}
\end{split}
\]

We notice that if $|\lambda_1|\not= |\lambda_2|$, a unique plane
is calibrated by $\alpha$. However if $|\lambda_1|=|\lambda_2|$,
then a two-dimensional family of planes in $G(2,V)$ are calibrated
by $\alpha$.

\end{proof}

\begin{lem}
A self-dual or anti-self-dual calibrating
form $\alpha$ can be written as
$\alpha(\cdot, \cdot)=<K(\cdot), \cdot>$
with $K\in O(4, V)$, the orthogonal group.
\end{lem}
\begin{proof}
If $\alpha $ is calibrating and
self-dual or anti-self-dual, then
$\max\{\lambda_1,\lambda_2\}=1 $ and
$\lambda_1 \lambda_2=\pm 1$, therefore $\lambda_1
=\pm 1$ and it is not hard to see that
$K$ is an isometry.

On the other hand, if $\alpha $ is induced
by an isometry $J$, then $\det J=\pm 1$,
therefore $\lambda_1\lambda_2=\pm 1$ and
$\alpha$ is self-dual or anti-self-dual.
\end{proof}

\begin{pro}\label{basis}
Let $(x, \mu)$ be an oriented two-plane in
$V$. Let $\alpha$ be a self-dual
calibrating form and $\beta$ be a anti
-self-dual calibrating form. Then there exists
an orthonormal basis $\{e_1,
e_2, e_3, e_4\}$ for $V$ with $\{e_1, e_2\}$
a basis for $x$ such that
$\mu(e_1, e_2)>0$ , $\nu(e_1, e_2, e_3, e_4)
>0$, $\alpha(e_A, e_B),
A, B=1\cdots 4 $ is of the form

\begin{equation}
\begin{pmatrix}0&\eta_1&\zeta_1&0\\
-\eta_1&0&0&-\zeta_1\\
-\zeta_1&0&0&\eta_1\\
0&\zeta_1&-\eta_1&0
\end{pmatrix}
\end{equation}
where $\eta_1=\alpha(e_1, e_2)$,
$\zeta_1^2+\eta_1^2=1$, and $\beta(e_A, e_B)$
is of the form.
\begin{equation}
\begin{pmatrix}0&\eta_2&\zeta_2&0\\
-\eta_2&0&0&\zeta_2\\
-\zeta_2&0&0&-\eta_2\\
0&-\zeta_2&\eta_2&0
\end{pmatrix}
\end{equation}
where $\eta_2=\beta(e_1, e_2)$ and
$\eta_2^2+\zeta_2^2=1$.
\end{pro}

\begin{proof}
Let $K, L$ be the elements in $End (V)$
corresponding to $\alpha$ and $\beta$.
If $\eta_1 \not=\pm 1$, we take any orthonormal basis
$\{e_1, e_2\}$ with $\mu(e_1, e_2)>0$. Notice
that $(K e_1)^T=\eta_1 e_2$.
Let
\[
\begin{split}
&e_3=\frac{1}{\sqrt{1-\eta_1^2}}
(K e_1-\eta_1 e_2)\\
&e_4=\frac{-1}{\sqrt{1-\eta_1^2}}
(K e_2+\eta_1 e_1)
\end{split}
\]

Therefore $\alpha_{A, B}$ is
of the required form. If $\eta_1=\pm 1$,
then any $\{e_1, e_2, e_3, e_4\}$ compatible
with $\mu $ and $\nu$ works.

It is not hard to check that $K$ and
$L$ as elements in $End(V)$ commute
and $KL $ is a self-adjoint operator.
Therefore we can rotate $\{e_1, e_2\}$ to
get a new basis so that
$<KL e_1, e_2>=0$.

This implies
\[<L e_1 ,e_4>
={\frac{-1}{\sqrt{1-\eta_1^2}}}<L e_1,
K e_2+\eta_1 e_1>=
{\frac{1}{\sqrt{1-\eta_1^2}}}
<KL e_1, e_2>=0
 \]
Likewise $<L e_2, e_3>=0$.
That $<L e_1, e_3>=-<L e_2, e_4>$
follows from the fact that $\beta $ is
anti-self-dual.

\end{proof}

Finally, we make a remark about $\alpha+\beta$.
In the above basis $\alpha+\beta$ is of the
form
\begin{equation}
\begin{pmatrix}0&\eta_1+\eta_2&\zeta_1+\zeta_2&0\\
-\eta_1-\eta_2&0&0&-\zeta_1+\zeta_2\\
-\zeta_1-\zeta_2&0&0&\eta_1-\eta_2\\
0&\zeta_1-\zeta_2&-\eta_1+\eta_2&0
\end{pmatrix}
\end{equation}

If the eigenvalues of $\sqrt{-1}(\alpha+\beta)$
are $\pm \lambda_1$ and $ \pm \lambda_2$, then
it is not hard to compute that
$\lambda_1^2\lambda_2^2=((\eta_1+\eta_2)
(\eta_1-\eta_2)-(\zeta_1+\zeta_2)
(-\zeta_1+\zeta_2))^2=0$ and $\lambda_1^2
+\lambda_2^2=(\eta_1+\eta_2)^2
+(\zeta_1+\zeta_2)^2
+(\eta_1-\eta_2)^2
+(\zeta_1-\zeta_2)^2=4$.  Therefore
$\frac{1}{2}(\alpha+\beta)$ is a calibrating
form and calibrates a unique two-plane.

\section{Surfaces in K\"ahler manifolds }

In the section, we assume $\baromega$
is a parallel self-dual calibrating two form
and $\baromega(X,Y)=<J(X),Y>$. $J$ is then
a parallel almost complex structure.
$M$ is therefore a K\"ahler manifold with
K\"ahler form $\baromega$.

We shall compute
 the equation of $\eta_t
=*_t\omega_t$ along the mean curvature flow.

The following Lemma is well-known.
\begin{lem}
Let $K(\cdot, \cdot)$ be the curvature operator
of $M$ and $Ric(\cdot, \cdot)$ be the Ricci
tensor of $M$. In terms of any orthonormal basis
$\{e_1, e_2, e_3, e_4\}$, the Ricci form is

\[Ric(JX,Y)={1\over 2} K(X,Y, e_A, J(e_A))\]
\end{lem}

\begin{proof}

This is seen by the following calculation

\[
\begin{split}
&K(JX, e_A, Y, e_A)\\
=&K(JX, e_A, J(Y), J(e_A))\\
=&-K(JX, JY, J(e_A), e_A)
-K(JX, J(e_A), e_A, J(Y))\\
=&K(X,Y, e_A, J(e_A))
-K(JX, J(e_A), Y, J(e_A))\\
\end{split}
\]
Now $ K(JX, e_A, Y, e_A)
=K(JX, J(e_A), Y, J(e_A))$ since $\{J(e_A)\}$ is
also an orthonormal basis.
\end{proof}

Let $F:\Sigma \rightarrow M$ be an isometric immersion. $\Sigma $
is equipped with a fixed orientation $d\mu$. By Proposition
\ref{basis}, for any point $p\in \Sigma$ it is possible to choose
an orthonormal basis $\{e_1, e_2, e_3, e_4\}$ for $T_p M$ such
that $d\mu(e_1, e_2)>0$ and $\baromega^2 (e_1, e_2, e_3, e_4)=
\baromega(e_1, e_2)\baromega(e_3, e_4)
-\baromega(e_1,e_3)\baromega(e_2, e_4) +\baromega(e_1,
e_4)\baromega(e_2,e_3)
>0$
and such that
$\baromega_{A, B}=\baromega(e_A, e_B),
A,B=1 \cdots 4$ is of the form.

\begin{equation}\label{momega}
\begin{pmatrix}0&\eta&\sqrt{1-\eta^2}&0\\
-\eta&0&0&-\sqrt{1-\eta^2}\\
-\sqrt{1-\eta^2}&0&0&\eta\\
0&\sqrt{1-\eta^2}&-\eta&0
\end{pmatrix}
\end{equation}
where $\eta=\omega(e_1, e_2)=*\omega$.

We first use this basis to calculate the curvature term in
Proposition (\ref{star}).

\begin{pro}\label{eta}
Let $Ric(\cdot, \cdot)$ be the Ricci tensor of
$M$. $\baromega$ a parallel
K\"ahler form.  Then $\eta=*_t\omega_t$ satisfies the following
equation

\begin{equation}\label{eta1}
\begin{split}
\frac{d}{dt}\eta=\Delta\eta+\eta[(h_{31k}
-h_{42k})^2
+(h_{32k}+h_{41k})^2]+(1-\eta^2) Ric(Je_1, e_2)
\end{split}
\end{equation}
\noindent
where $\{e_1, e_2, e_3,e_4\}$
is any orthonormal basis for $T_p M$ such that
$\{e_1, e_2\}$ forms an orthonormal basis
for $T\Sigma_t$, $d\mu(e_1, e_2)>0$ and $\baromega^2
(e_1, e_2, e_3, e_4)>0$. $A(e_i,e_j)=h_{3ij}e_3+h_{4ij}e_4$
is the second fundamental form.
\end{pro}

\begin{rem}
Notice that the term $(h_{31k}
-h_{42k})^2
+(h_{32k}+h_{41k})^2$ depends only on the
orientation of $\{e_1, e_2, e_3, e_4\}$ but
not on the particular orthonormal basis
we choose.
\end{rem}
\begin{proof}

First we show
\begin{equation}\label{ric}
\baromega((K(e_k, e_1)e_k)^N,e_2)
-\baromega((K(e_k, e_2)e_k)^N,e_1)
=(1-\eta^2)Ric(Je_1, e_2)
\end{equation}

By definition,
\[
\begin{split}
&\baromega((K(e_k, e_1)e_k)^N,e_2)
-\baromega((K(e_k, e_2)e_k)^N,e_1)\\
=-&<(Je_2)^N,K(e_k, e_1)e_k>
+<(Je_1)^N,K(e_k, e_2)e_k>
\end{split}
\]
Therefore when $\eta=\pm 1$, equation (\ref{ric})
is obvious, therefore we may assume $\eta
\not= \pm 1$ and apply the basis in equation
(\ref{momega}) and get

\[
\sqrt{1-\eta^2}(K(e_k, e_2,e_k,e_3)
+K(e_k, e_1, e_k, e_4))
=\sqrt{1-\eta^2}(K(e_1, e_2,e_1,e_3)
+K(e_2, e_1, e_2, e_4))
\]

By the previous lemma,

\[
\begin{split}
& Ric(JX, Y)\\
&={1\over 2} J_{AB}K(X,Y, e_A, e_B)\\
&=\eta  (K(X,Y,e_1,e_2)+
K(X,Y, e_3,e_4))
+\sqrt{1-\eta^2} (K(X,Y,e_1,e_3)-K
(X,Y, e_2,e_4))
\end{split}
\]

Since $J$ is parallel and isometry, the curvature tensor
is $J$ invariant, therefore we have

\[K(X,Y, e_1, e_2)=K(X,Y, J(e_1)
,J(e_2))\]

Use (\ref{momega}) again, this is the same as

\[(1-\eta^2)(K(X,Y, e_1, e_2)
+K(X,Y, e_3, e_4))
=\eta\sqrt{1-\eta^2} (K(X,Y, e_1, e_3)
-K(X,Y, e_2, e_4))
\]

Therefore
\[
\begin{split}
Ric(JX,Y)&=\frac{1}{\sqrt{1-\eta^2}}
(K(X,Y, e_1,e_3)-K(X,Y, e_2, e_4))\\
\end{split}
\]

Equation (\ref{ric}) now follows by
substituting $X=e_1, Y=e_2$.

We use the basis in equation (\ref{momega})
to calculate the rest terms  in Proposition
\ref{star}.

\[
\begin{split}
&\eta|A|^2-2\baromega(A(e_k,e_1), A(e_k, e_2))\\
=&\eta(|A|^2-2h_{31k}h_{42k}+2h_{41k}
h_{32k})
\end{split}
\]
Equation (\ref{eta1}) follows by completing squares.

\end{proof}

\begin{rem}
When $M$ is a K\"ahler manifold with
K\"ahler form $\baromega$ and almost
complex structure $J$.  The second fundamental
form of a holomorphic submanifold has the
symmetry $h_{31k}=h_{42k}, h_{41k}
=-h_{32k}$.
\end{rem}

\section{Asymptotics of Singularities}

In this section, we study the asymptotic behavior of singularities
of the mean curvature flow. In particular, we show that no type I
singularity will occur in the mean curvature flow of symplectic
surface in a four-dimensional K\"ahler-Einstein manifold.
Techniques involved are blow-up analysis and monotonicity formula
of backward heat kernel.

The following lemma
says singularity forms only when the second
fundamental form blows up.

\begin{lem}
Given any mean curvature flow $F:\Sigma
\times [0,t_0)
\rightarrow M $, suppose
$\sup_{t\in[0,t_0)}\sup_{x\in \Sigma}
 |A|(x,t)$ is bounded where
$|A|(x,t)$ is the norm of the second
fundamental form for $F_t(\Sigma)$ at $F_t(x)$. Then $F$
can be extended to $\Sigma\times [0,\bar{t}_0)$
for some $\bar{t}_0> t_0$.
\end{lem}

\begin{proof}
It can be shown that all higher covariant
derivatives of the second fundamental form are
uniformly bounded. For the detail see \cite
{hu1} for
the hypersurface case.
\end{proof}

Since the study of singularities is local,
it is more convenient to adopt an unparametrized definition
of mean curvature flow introduced in
\cite{w3}. Let $M$ be an $m-$dimensional Riemannian
manifold of bounded geometry. An
immersed smooth submanifold $\mathcal{S}\subset M\times \R$
is a smooth flow if the function $\tau
:M\times \R\rightarrow \R$, $\tau(y,t)=t$ has no
critical points in $\mathcal{S}$.
 $\mathcal{S}_t=\mathcal{S}\cap M\times\{t\}$ is
 called the $t$-slice of $\mathcal{S}$.
At each point $(y,t)\in \mathcal{S}$, the normal velocity $v(y,t)$
is the unique vector that satisfies $v$ is normal to
$\mathcal{S}_t$ and $v+\frac{\partial}{\partial t}$ is tangent to
$\mathcal{S}$. $H(y,t)$ is the mean curvature vector of
$\mathcal{S}_t$ at $(y,t)$. We allow $M$ and $\mathcal{S}_t$ to
have boundary. In fact, all unparametrized flow considered in this
article is of the form $\cup_{t\in [0, t_0)} (F_t(\Sigma)\cap
B)\times \{t\}$, where $F$ is a parameterized mean curvature flow
of a compact manifold $\Sigma$ without boundary and $B$ is an
neighborhood of $y$ in a complete Riemannian manifold. Therefore
$\partial \mathcal{S}_t \subset \partial M$.

\begin{dfn}
A smooth flow $\mathcal{S}$ is called a (unparametrized)
mean curvature flow if
\[ v(y,t)=H(y,t)\]
at each point of $(y,t)\in \mathcal{S}$.
\end{dfn}

Let $(y_0, t_0)$ be an interior point in $M\times \R $. When $M$
is the Euclidean space, in \cite{hu2} Huisken introduces the
backward heat kernel to study the asymptotic behavior near
singular points. Recall the (n-dimensional) backward heat kernel
$\rho_{y_0, t_0}$ at $(y_0, t_0)$.

\begin{equation}
\rho_{y_0, t_0}(y,t)=\frac{1}{(4\pi(t_0-t))^{n\over 2}}
\exp (\frac{-|y-y_0|^2}{4(t_0-t)})
\end{equation}
The monotonicity formula of Huisken asserts
for $t<t_0$

\[\frac{d}{dt}\int \rho_{y_0, t_0} d\mu_t
\leq 0\]

 For general Riemannian manifold $M$, following
\cite{w2}, we isometrically embed $M$ into
$\R^N$. The mean curvature flow of $\Sigma$ in
$M$ now reads.

\[\frac{d}{dt} F=H=\overline{H}+E\]
where $F$ is the coordinate function in $\R^N$,
$H$ is the mean curvature vector of $\Sigma$ in
$M$, $\overline{H}$ is the mean curvature vector
of $\Sigma$ in $\R^N$, and

\[E=\sum_i \overline{A}(e_i, e_i)\]

Here $\overline{A}$ denotes the second fundamental form of $M$ in
$\R^N$ and $\{e_i\}$ is an orthonormal basis for $T\Sigma_t$.

In the general case $\int \rho_{y_0,t_0} d\mu_t$
is no longer decreasing, however
the following is still true.

\begin{pro}\label{mono}
Let $\mathcal{S}\subset M\times \R$ be
a mean curvature flow such that $\partial
\mathcal{S}_t \subset \partial M$. We
fix an isometric embedding $M\hookrightarrow
\R^N$ and let $\rho_{y_0,t_0}$ be
the (n-dimensional) backward heat kernel
at $(y_0, t_0)$. Then the limit

\[\lim_{t\rightarrow t_0} \int \rho_{y_0, t_0}
d\mu_t\]
exists, where $d\mu_t$ is the Radon measure
associated with $\mathcal{S}_t \subset M$.
\end{pro}

\begin{proof}
See Proposition 11 in \cite{w2}.
\end{proof}

The limit is called the Gaussian density
of $\mathcal{S}$ at $(y_0, t_0)$ in \cite{w3}.
The Gaussian density can be used to detect
singularities of mean curvature flow.
The following theorem of White in \cite{w3}
is a parabolic analogue of Allard's regularity theorem.

\begin{thm}
There is an $\epsilon >0$ such that whenever

\[\lim_{t\rightarrow t_0}\int \,\rho_{y_0, t_0}
d\mu_t <1+\epsilon\], it can concluded that $(y_0, t_0)$ is a
regular point of $\mathcal{S}$.
\end{thm}

A regular point is a point where the second
fundamental form is locally bounded in H\"older
norm.

To study singularity, we consider the parabolic blow-up near a
possible singular point. Let $F:\Sigma \times [0,t_0) \rightarrow
M\hookrightarrow \R^N $ be a parameterized mean curvature flow.
Let $B$ be a ball about $y_0$ of radius $r$ in $\R^N$.
 Take $\mathcal{S}=\cup_{t\in [0, t_0)}(F_t(\Sigma)\cap B)\times \{t\}$,
  then $\mathcal{S}$
 is an unparametrized mean curvature flow
 in $B$.

For any $\lambda >1$, the parabolic dilation $D_\lambda$
at $(y_0, t_0)$ is
defined by

\begin{equation}
\begin{split}
D_\lambda :\,\R^N \times[0, t_0) &\rightarrow
\R^N \times [-\lambda^2 t_0, 0)\\
(y,t)&\rightarrow (\lambda(y-y_0), \lambda^2
(t-t_0))
\end{split}
\end{equation}

For any $s$, $-\lambda ^2 t_0\leq s<0$, the
two slices
$\mathcal{S}^\lambda_s$ and $\mathcal{S}
_{t_0+{s\over \lambda^2}}$ can be identified
and $d\mu_s^\lambda=\lambda^n d\mu_t$.

It is not hard to check that if we denote $F_s^\lambda (x)=\lambda
(F_t(x)-y_0)$ for $s=\lambda^2(t-t_0)$, then

\[\rho_{0,0}(F_s^\lambda(x), s)=\rho_{0,0}(\lambda(F_t(x)-y_0),
\lambda^2(t-t_0))=\frac{1}{\lambda^n} \rho_{y_0,t_0}(F_t(x),t)\]

Therefore
\[\int \rho_{y_0,t_0} d\mu_t
=\int \rho_{0,0} d\mu_s^\lambda
\]
is invariant under the parabolic dilation.

 The singularity
of $\mathcal{S}$ near $(y_0, t_0)$ is reflected in the asymptotic
behavior of $\mathcal{S}^\lambda$ as $\lambda \rightarrow \infty$.

Take any sequence $\lambda_i\rightarrow
\infty$, it can be proved as in \cite{il1} and
\cite{w2} that
a subsequence of $\mathcal{S}^{\lambda_i}$
converges to a Brakke flow
$\mathcal{S}^\infty\subset \R^N \times
(-\infty, 0)$. $\mathcal{S}^\infty$ is
called a tangent flow of $\mathcal{S}$ at
$(y_0, t_0)$.

Now we state and prove the main proposition in
this section.

\begin{pro}\label{stat}
If $F:\Sigma\times [0,t_0)\rightarrow M
\hookrightarrow \R^N $
is a mean curvature flow of an
 orientable surface
in a (real) four dimensional K\"ahler
manifold $M$. Assume the second fundamental
form of $M \hookrightarrow \R^N$ is bounded.
Let $\omega(\cdot, \cdot)
=<J(\cdot), \cdot>$ be a K\"ahler form
on $M$.  If
there exist $\delta, C >0$ such that
$\eta_t=* \omega_t>\delta$ on $F_t(\Sigma)$ for
$t\in [0 ,t_0)$ and such that $|A|^2\leq \frac{C}
{t_0-t}$, then $F$
can be extended to $\Sigma\times [0,\bar{t}_0)$
for some $\bar{t}_0> t_0$.
.
\end{pro}

\begin{proof}
Let $y_0\in M$, we shall consider the blow
up of the mean curvature flow at $(y_0, t_0)$.
Let $B$ be a ball of radius $r$ about
$y_0$ in $\R^N$ and $\psi$ be a cut-off
function supported in $B$ so that $\psi\equiv 1$
in the ball of radius $\frac{r}{2}$ about
$y_0$. We assume
\[|\overline{\nabla}\psi|+|\overline{\nabla}
\overline{\nabla}\psi|\leq C\]
where $\overline{\nabla}$ is the covariant
derivative on $\R^N$.
Recall the equation for $\eta$
is
\[
\begin{split}
\frac{d}{dt}\eta&=\Delta\eta+\eta[(h_{31k}
-h_{42k})^2
+(h_{32k}+h_{41k})^2]+(1-\eta^2) Ric(Je_1, e_2)\\
\end{split}
\]

The backward heat kernel $\rho_{y_0,t_0}$ satisfies the following
parabolic equation along the mean curvature flow. Notice that
$\nabla$ and $\Delta$ are the covariant derivative and the Laplace
operator on $\Sigma_t$ respectively.

\begin{equation}\label{bh}
\begin{split}
\frac{d}{dt}\rho_{y_0,t_0}&=-
\Delta\rho_{y_0,t_0}
-\rho_{y_0,t_0}( \frac{|F^\perp|^2}{4(t_0-t)^2}
+\frac{F^\perp \cdot \overline{H}}{t_0-t}
+\frac{F^\perp \cdot E}{2(t_0-t)})
\end{split}
\end{equation}
where $F^\perp$ is the component of $F\in T\R^N$ in $T\R^N \slash
T\Sigma_t$. This equation for mean curvature flow in a Euclidean
space is essentially derived by Huisken \cite{hu2} and in a
general ambient manifold by White \cite{w2}. It is derived in the
next paragraph for completeness. Recall that
\[ \frac{d}{dt} F(x,t)=H=\overline{H}+E\]
where $H\in TM\slash T\Sigma$ is the mean curvature vector of
$\Sigma_t$ in $M$ and $\overline{H}\in T\R^N\slash T\Sigma$ is the
mean curvature vector of $\Sigma_t$ in $\R^N$.

We may assume $y_0$ is the origin and then

\[
\rho_{y_0, t_0}(F(x,t),t)=\frac{1}{(4\pi(t_0-t))^{n\over 2}} \exp
(\frac{-|F(x,t)|^2}{4(t_0-t)})
\]

Abbreviate $\rho_{y_0, t_0}(F(x,t),t)$ by $\rho$, it is not hard
to see

\begin{equation}\label{ddtrho}
\frac{d}{dt}\rho
=\rho[\frac{n}{2(t_0-t)}-\frac{|F(x,t)|^2}{4(t_0-t)^2}-\frac{F\cdot
H}{2(t_0-t)}]
\end{equation}

We shall compute $\sum_i (\barbla_{e_i} \barbla \rho)\cdot e_i$ in
two different ways, where $\barbla$ denotes the covariant
derivative in $\R^N$ and $\{e_i\}$ is an orthonormal basis for
$T\Sigma$.

\[\sum_i (\barbla_{e_i} \barbla \rho)\cdot e_i=\sum_i \barbla_{e_i} (\nabla \rho
+(\barbla \rho)^{T\R^N\slash T\Sigma})\cdot e_i=\Delta \rho
-\barbla \rho\cdot \overline{H}
\]
$\barbla\rho=-\frac{\rho}{2(t_0-t)} F$, thus
\[\sum_i (\barbla_{e_i} \barbla \rho)\cdot e_i=\Delta \rho
+\frac{1}{2(t_0-t)} F\cdot\overline{H}
\]

On the other hand,
\[
\begin{split}
&\sum_i (\barbla_{e_i} \barbla \rho)\cdot e_i\\
&=\sum_i \barbla_{e_i} (-\frac{\rho}{2(t_0-t)} F)\cdot e_i\\
&= -\frac{1}{2(t_0-t)}(\nabla \rho\cdot F +\rho
\sum_i\barbla_{e_i}
F\cdot e_i)\\
&=-\frac{1}{2(t_0-t)}(-\frac{\rho}{2(t_0-t)}F^{T\Sigma}\cdot F +n
\rho)\\
&=\rho(\frac{1}{4(t_0-t)^2}|F^{T\Sigma}|^2-\frac{n}{2(t_0-t)})
\end{split}
\]

Compare these two equalities, we get

\begin{equation}\label{laprho}
\Delta \rho=\rho[-\frac{1}{2(t_0-t)}F\cdot \overline{H}
+\frac{1}{4(t_0-t)^2}|F^{T\Sigma}|^2-\frac{n}{2(t_0-t)}]
\end{equation}

Now add equations (\ref{ddtrho}) and (\ref{laprho}), we get

\[
\begin{split}
&\frac{d}{dt}\rho +\Delta
\rho\\
&=\frac{\rho}{4(t_0-t)^2}(|F^{T\Sigma}|^2-|F|^2)-\frac{\rho}{2(t_0-t)}
(F\cdot \overline{H}+F\cdot H)\\
&=-\frac{\rho}{4(t_0-t)^2}|F^\perp|^2-\frac{\rho}{2(t_0-t)}
(F^\perp\cdot \overline{H}+F^\perp \cdot H)
\end{split}
\]
where $F^\perp=(F)^{T\R^N\slash T\Sigma}$. Recall that
$H=\overline{H}+E$ and we get equation (\ref{bh}).

 The minus sign in front of the Laplacian in equation (\ref{bh})
indicates the fact that $\rho$ satisfies the backward heat
equation. The following inequality is particularly useful when
deal with backward heat kernels.

\begin{equation}\label{iden}
 g(-\Delta\rho)+(\Delta g)\rho
=-\text{div}(\nabla \rho\, g)
+\text{div}(\rho \nabla g)
\end{equation}

The volume form $d\mu_t$ of $\Sigma_t$ satisfies the equation

\[
\frac{d}{dt}d\mu_t
=-|H|^2 d\mu_t
=-\overline{H}\cdot (\overline{H}+E)
d\mu_t
\]

Therefore,
\[
\begin{split}
&\frac{d}{dt}\int \psi(1-\eta)
\rho_{y_0,t_0}\, d\mu_t\\
&= \int [\frac{d}{dt}\psi(1-\eta)]
\rho_{y_0,t_0}\, d\mu_t +\int \psi(1-\eta)
[\frac{d}{dt}\rho_{y_0,t_0}]\, d\mu_t
-\int \psi(1-\eta)
\rho_{y_0,t_0}\,\overline{H}\cdot (\overline{H}+E) d\mu_t
\end{split}
\]

Plug the equation (\ref{bh}) for $\frac{d}{dt}\rho_{y_0,t_0}$ ,
use the identity (\ref{iden}) with $g=\psi(1-\eta)$, and complete
square we get

\begin{equation}
\begin{split}
&\frac{d}{dt}\int \psi(1-\eta)
\rho_{y_0,t_0}\, d\mu_t\\
=&\int[\frac{d}{dt} (\psi(1-\eta))-\Delta
(\psi(1-\eta))]\rho_{y_0,t_0} \,d\mu_t \\
&-\int \psi(1-\eta)\rho_{y_0,t_0}[|\overline{H}
+\frac{1}{2(t_0-t)}F^\perp|^2
+(\overline{H}
+\frac{1}{2(t_0-t)}F^\perp)\cdot E ]
d\mu_t
\end{split}
\end{equation}

Now
\[\frac{d}{dt} (\psi(1-\eta))-\Delta
(\psi(1-\eta)) =\psi(-\frac{d}{dt}\eta+\Delta \eta) +(\barbla
\psi\cdot H)(1-\eta) +2\nabla\psi \cdot \nabla \eta -\Delta \psi
(1-\eta)
\]
where we use $\frac{d}{dt}\psi=\barbla\psi \cdot H$.

Integration by parts,
\[
\int [2\nabla\psi \cdot \nabla \eta
-\Delta \psi (1-\eta)]\rho_{y_0,t_0} d\mu_t
=\int[\nabla \psi\cdot \nabla \eta \,\rho_{y_0,t_0}
+\nabla \psi \cdot \nabla \rho_{y_0, t_0} \,(1-\eta)]
d\mu_t
\]

Therefore, we have
\[
\begin{split}
&\frac{d}{dt}\int \psi (1-\eta) \rho_{y_0,t_0}\,
d\mu_{t}\\
=&-\int \psi \eta \rho_{y_0, t_0}
[(h_{31k}-h_{42k})^2+(h_{32k}+h_{41k})^2]
\,d\mu_t\\
&-\int\psi(1- \eta)\rho_{y_0, t_0} |\overline{H}+{1\over 2(t_0
-t)} F^\perp +\frac{E}{2}|^2 \,d\mu_t +
\int\psi (1-\eta)\rho_{y_0,t_0}\frac{|E|^2}{4}\, d\mu_t\\
&+\int [(\barbla \psi \cdot H)(1-\eta) \rho_{y_0,t_0}+\nabla
\psi\cdot \nabla\eta \,\rho_{y_0,t_0} +\nabla \psi \cdot
\nabla\rho_{y_0,t_0}\, (1-\eta)]d\mu_t
\end{split}
\]

Since $|E|$ and $\int \rho_{y_0,t_0}d\mu_t $ are both bounded,

\[
\begin{split}
&\frac{d}{dt}\int \psi (1-\eta) \rho_{y_0,t_0}\,
d\mu_{t}\\
\leq & C-\int \psi \eta \rho_{y_0, t_0}
[(h_{31k}-h_{42k})^2+(h_{32k}+h_{41k})^2]
\,d\mu_t\\
&+\int [(\barbla \psi \cdot H)(1-\eta) \rho_{y_0,t_0}+\nabla
\psi\cdot \nabla\eta \,\rho_{y_0,t_0} +\nabla \psi \cdot
\nabla\rho_{y_0,t_0}\, (1-\eta)]d\mu_t
\end{split}
\]

The last term is also bounded by the following
computation.
\[\int \nabla \psi \cdot \nabla \rho_{y_0, t_0} \,(1-\eta)\,d\mu_t
\leq C \int_{B\backslash B_{{1\over 2}r} (y_0)}|\nabla\rho_{y_0,
t_0}| d\mu_t
\]

Since $\nabla\rho_{y_0,t_0}=-\rho_{y_0, t_0} \frac{\nabla
\,|F-y_0|^2}{4(t_0-t)}$ and $|\nabla \,|F-y_0|^2|\leq  |\barbla\,
|F-y_0|^2|| \leq 2|F-y_0|$, we have

\[\int \nabla \psi \cdot \nabla \rho_{y_0, t_0} \,(1-\eta)\,d\mu_t
\leq C \int_{B\backslash
B_{{1\over 2}r} (y_0) }
\frac{1}{(t_0-t)^{{n\over2}+1}}
\exp (\frac{-{{1\over 4}r}^2}{4(t_0-t)}) d\mu_t
\]
The last expression approaches zero as
$t\rightarrow t_0$.

Therefore

\[
\begin{split}
&\frac{d}{dt}\int \psi (1-\eta) \rho_{y_0,t_0}\,
d\mu_{t}\\
\leq & C-\int \psi \eta \rho_{y_0, t_0}
[(h_{31k}-h_{42k})^2+(h_{32k}+h_{41k})^2]
\,d\mu_t\\
&+\int (\barbla \psi \cdot H)(1-\eta) \rho_{y_0,t_0} d\mu_t+ \int
\nabla \psi\cdot \nabla\eta \,\rho_{y_0,t_0}d\mu_t
\end{split}
\]

The term
$\int \nabla \psi\cdot \nabla\eta
 \,\rho_{y_0, t_0} \,d\mu_t$ can be written in  the
 following

\[ \int (\frac{\nabla \psi}{\sqrt{\psi}}
\sqrt{\rho_{y_0, t_0}})\cdot (\nabla\eta\, \sqrt{\psi}
 \,\sqrt{\rho_{y_0, t_0}}) \,d\mu_t
 \leq \frac{1}{4\epsilon^2}\int
  \frac{|\barbla \psi|^2}{\psi}
{\rho_{y_0, t_0}}\, d\mu_t +\epsilon^2 \int |\nabla \eta|^2 \psi
\rho_{y_0, t_0} d\mu_t
\] where we use $ |\nabla \psi|^2\leq |\barbla
\psi|^2$.

In a normal coordinate system, we compute
$\nabla\eta$, again use the basis in
equation (\ref{momega}).

\begin{equation}\label{deta}
\begin{split}
\partial_k \eta
=&\partial_k(\frac{\omega(\partial_1,
\partial_2)}{\sqrt{\det g}})\\
=&\pk(\omega(\partial_1, \partial_2))\\
=&\omega(A(\pk, \partial_1), \partial_2)
+\omega(\partial_1, A(\pk, \partial_2))\\
=&h_{\alpha 1k}{\omega}_{\alpha 2}
+h_{\alpha 2k}{\omega}_{1\alpha}\\
=&\sqrt{1-\eta^2}(h_{41k}+h_{32k})
\end{split}
\end{equation}

Therefore $|\nabla\eta|^2\leq (1-\eta^2) (h_{41k}+h_{32k})^2$ and
thus
\[\int \nabla \psi\cdot \nabla\eta
 \,\rho_{y_0, t_0} \,d\mu_t \leq \frac{1}{4\epsilon^2}\int
  \frac{|\barbla \psi|^2}{\psi}
{\rho_{y_0, t_0}}\, d\mu_t +\epsilon^2 \int
(1-\eta^2)(h_{41k}+h_{32k})^2 \psi \rho_{y_0, t_0} d\mu_t
 \]
Likewise since $|H|^2\leq
2[(h_{31k}-h_{42k})^2+(h_{32k}+h_{41k})^2]$, we have

\[
\begin{split}
&\int (\barbla \psi \cdot H)(1-\eta) \rho_{y_0,t_0} d\mu_t \\
&\leq \frac{1}{4\epsilon^2}\int \frac{|\barbla \psi|^2}{\psi}
{\rho_{y_0, t_0}}\, d\mu_t +2\epsilon^2 \int
(1-\eta)^2[(h_{31k}-h_{42k})^2+(h_{41k}+h_{32k})^2] \psi
\rho_{y_0, t_0} d\mu_t
\end{split}
\]
 Since $\psi$ is of compact support, by Lemma 6.6
in \cite{il2}, $\frac{|\barbla \psi|^2}{\psi}\leq 2 \max |\barbla
\barbla \psi|$ is bounded .


Since $\eta>\delta$, we can choose $\epsilon $ small enough so
that

\begin{equation}\label{5.5}
\begin{split}
&\frac{d}{dt}\int \psi (1-\eta) \rho_{y_0,t_0}\,
d\mu_{t}\\
\leq &\, C-C_\delta \int  \psi \rho_{y_0, t_0}
[(h_{31k}-h_{42k})^2+(h_{32k}+h_{41k})^2] \,d\mu_t
\end{split}
\end{equation}
where $C_\delta$ is a constant that depends on $\delta$.

From this we see that
$\lim_{t\rightarrow t_0}
\int \psi (1-\eta)\rho_{y_0, t_0} d\mu_t$
exists.

For $\lambda>1$, let's study the flow $\mathcal{S}^\lambda \subset
\R^N \times [-\lambda^2 t_0, 0)$. Let $\rho^\lambda_{0,0,}(y,s)$
be the backward heat kernel  at $(0,0)$ and
$\psi^\lambda(F_s^\lambda(x))=\psi (F_t(x))$. Recall that
$t=t_0+\frac{s}{\lambda^2}$, thus
\[
\begin{split}
&\frac{d}{ds}\int \psi^\lambda
 (1-\eta^\lambda) \rho^\lambda_{0,0}\,
d\mu_{s}^{\lambda}\\
=&\frac{1}{\lambda^2}\frac{d}{dt}\int \psi (1-\eta) \rho_{y_0,t_0}\,
d\mu_{t}\\
\leq &\, \frac{C}{\lambda^2}-\frac{C_\delta} {\lambda^2} \int
\psi \rho_{y_0, t_0} [(h_{31k}-h_{42k})^2+(h_{32k}+h_{41k})^2]
\,d\mu_t\\
\end{split}
\]
We notice that $\eta $ is a scaling invariant
quantity therefore $\eta^\lambda=\eta$.
It is not hard to check that
\[
\frac{1} {\lambda^2} \int  \psi \rho_{y_0, t_0}
[(h_{31k}-h_{42k})^2+(h_{32k}+h_{41k})^2]\,d\mu_t =\int
\psi^\lambda \rho^\lambda_{0,0}
[(h^\lambda_{31k}-h^\lambda_{42k})^2+(h^\lambda_{32k}+h^\lambda_{41k})^2]
\,d\mu_{s}^{\lambda}
\]

This is because $\rho_{y_0, t_0} d\mu_t$
is invariant under the parabolic
scaling and the norm of second fundamental
form scales like the inverse of
the distance.

Therefore

\[
\begin{split}
&\frac{d}{ds}\int \psi^\lambda
 (1-\eta^\lambda) \rho^\lambda_{0,0}\,
d\mu_{s}^{\lambda}\\
\leq &\frac{C}{\lambda^2}-C_\delta \int \psi^\lambda
\rho^\lambda_{0,0}
[(h^\lambda_{31k}-h^\lambda_{42k})^2+(h^\lambda_{32k}+h^\lambda_{41k})^2]\,
d\mu_s^\lambda
\end{split}
\]

Compare with equation (\ref{5.5})
and we see this reflect the correct
scaling for the parabolic blow-up.

Take any $\tau>0$ and integrate from
$-1-\tau$ to $-1$.

\begin{equation}\label{delta}
\begin{split}
&C_\delta\int_{-1-\tau}^{-1} \int \psi^\lambda \rho^\lambda_{0,0}
[(h^\lambda_{31k}-h^\lambda_{42k})^2+
(h^\lambda_{32k}+h^\lambda_{41k})^2]
\,d\mu_{s}^{\lambda}ds\\
&\leq \int \psi^\lambda
 (1-\eta^\lambda) \rho^\lambda_{0,0}\,
d\mu_{-1}^{\lambda}-\int \psi^\lambda
 (1-\eta^\lambda) \rho^\lambda_{0,0}\,
d\mu_{-1-\tau}^{\lambda}
+\frac{C}{\lambda^2}
\end{split}
\end{equation}

Notice that

\[\int \psi^\lambda (1-\eta^\lambda)
\rho_{0,0}^\lambda d\mu_{s}
^{\lambda}
=\int \psi (1-\eta)\rho_{y_0, t_0}
d\mu_{t_0+\frac{s}{\lambda^2}}
\]

This equality means the quantity $\int \psi (1-\eta)\rho_{y_0,
t_0} d\mu_{t}$ is invariant under parabolic scaling. This fact is
extremely important in applying the Monotonicity formula. Recall
the natural Monotonicity formula for the volume

\[\frac{d}{dt}\int d\mu_t
=-\int|H|^2 d\mu_t\]

But $\int d\mu_s^\lambda
=\lambda^2\int d\mu_t$ is not
scaling invariant. This deteriorates the
usefulness of the formula in the blow-up
analysis.

Now the right hand
side in equation (\ref{delta}) tends to
zero as $\lambda \rightarrow \infty$. For any sequence $\lambda_i\rightarrow \infty$,
we can choose $s_i\rightarrow -1$ such that
\[
\begin{split}
& \int \psi^{\lambda_i}
\rho^{\lambda_i}_{0,0}
[(h^{\lambda_i}_{31k}-h^{\lambda_i}_{42k})^2+
(h^{\lambda_i}_{32k}+h^{\lambda_i}_{41k})^2]
\,d\mu_{s_i}^{\lambda_i}\rightarrow 0
\end{split}
\]
as $i\rightarrow \infty$.

It is not hard to compute that

\[|A|^2(\mathcal{S}_{s}^{\lambda_i})
=\frac{1}{\lambda_i^2} |A|^2(\mathcal{S}_{t_0
+\frac{s}{\lambda_i^2}})
=(\frac{-1}{s})(t_0-t_i)|A|^2(\mathcal{S}_{t_i})\]

The assumption implies each $\Sigma_{s} ^{\lambda_i}$ has
uniformly bounded second fundamental form. By the same method used
in \cite{hu2}, any higher covariant derivatives of the second
fundamental form of $\mathcal{S}_{s}^ {\lambda_i}$ is bounded.
Therefore the convergence
$\mathcal{S}_{s_i}^{\lambda_i}\rightarrow \mathcal{S}_{-1}^\infty$
is smooth.

We may assume each $\mathcal{S}_t$ is connected by taking
connected components. Therefore we have
$(h_{31k}-h_{42k})^2+(h_{32k} +h_{41k})^2=0$ for
$\mathcal{S}_{-1}^\infty$. This implies $\nabla {\eta}=0$ and
$H=0$. Applying the same argument to the monotonicity formula for
$\int\psi \rho_{y_0,t_0}d\mu_t$ gives $H+\frac{1}{2}F^{\perp}=0$
for $\mathcal{S}_{-1}^\infty$. To sum up, we get $F^\perp=0$ and
$\nabla \eta=0$ for $\mathcal{S}^\infty_{-1}$. The first
 condition implies $\mathcal{S}_{-1}^\infty$ is a plane with multiplicity
 one. On the other hand,

\[
\lim_{t_i\rightarrow t_0}\int \rho_{y_0,t_0} d\mu_{t_i}
=\lim_{i\rightarrow \infty} \int
\frac{1}{4\pi(-s_i)}\exp({\frac{-|F^{\lambda_i}_{s_i}|^2}{4(-s_i)}})
\,d\mu^{\lambda_i}_{s_i} =\int
\frac{1}{4\pi}\exp({\frac{-|F^{\infty}_{-1}|^2}{4}})
\,d\mu^{\infty}_{-1}
\]
where $\lambda_i =\sqrt{\frac{-s_i}{t_0-t_i}}$.

The last Gaussian integral for a plane can
be calculated directly and is equal to $1$.  By White's theorem,
$(y_0, t_0) $ is a regular point.
\end{proof}

We recall the following definition of type I
singularities for the mean curvature flow.

\begin{dfn}
A singularity at $t_0$ is called type I if
there exists a $C$ such that $|A|^2\leq
\frac{C}{t_0-t}$.
\end{dfn}

Recall a K\"ahler manifold $M$ is
called  K\"ahler-Einstein if $Ric=c g$ for some constant $c$.
In this case, the scalar curvature $s$ of $M$ is
a constant and $c=\frac{s}{4}$.
A immersion $F:\Sigma\rightarrow M$ is symplectic
if $\eta>0$.

\vskip 10pt
\noindent {\bf Theorem A}
\hskip 3pt{\it Let $M$ be a four-dimensional
K\"ahler-Einstein manifold, then a symplectic
surface remains symplectic along the
mean curvature flow and the flow does not
develop any type I singularities.}
\vskip 10pt

\begin{proof}
Since $Ric (J\cdot, \cdot)
=\omega (\cdot, \cdot)$ by Proposition
\ref{eta}, the equation of $\eta=*\omega$ now becomes

\begin{equation}\label{eta2}
\frac{d}{dt}\eta =\Delta\eta
+\eta [(h_{31k}-h_{42k})^2
+(h_{32k}+h_{41k})^2+
c(1-\eta^2)]
\end{equation}

The first assertion follows from maximum principle for parabolic
equations. Actually, when $c\geq 0$, i.e. the nonnegative scalar
curvature case, the function $\min_{\Sigma}\,\eta_t$ is a
non-decreasing function of $t$.  In any case, by comparison
theorem for parabolic equations, $\eta$ has a positive lower bound
at any finite time and Proposition \ref{stat} is applicable.
\end{proof}

\begin{rem}
The same argument can be used to prove
there is no type I singularity for
the mean curvature flow of an almost
calibrated Lagrangian
submanifolds in a Calabi-Yau manifold
$M$.
Here the almost calibrated condition is
$*\Omega>0$ where $\Omega $ is the real
part of the canonical form on $M$.
In fact, $*\Omega$ satisfies

\[\frac{d}{dt}*\Omega =\Delta *\Omega
+|H|^2 *\Omega\]

A smooth blow-up limit satisfies $H+\frac{1}{2} F^\perp=0$ and
$H=0$ and is thus a linear subspace.
\end{rem}

\section{ Long time existence and convergence}

In this section, we study the problem of
long time existence. The main result is the following.

\begin{pro}\label{long}
Let $M$ be an oriented four-dimensional compact manifold. Let
$\omega'$ and $\omega''$ be two parallel calibrating form such
that $\omega'$ is self-dual and $\omega''$ is anti-self-dual . Let
$(\Sigma_0, d\mu)\hookrightarrow M $ be an compact surface with
orientation $d\mu$. Let $F:\Sigma \times [0,t_0)\rightarrow M $ be
the mean curvature flow of $\Sigma_0$ such that there exist a $
\delta>0$ with $*\omega'>\delta$ and $*\omega''>\delta$ on
$F_t(\Sigma)$ for $0\leq t<t_0$.  Then $F$ can be extended
smoothly to $\Sigma\times [0, \bar{t}_0)$ for some
$\bar{t}_0>t_0$.
\end{pro}

\begin{proof}
The assumption implies $(\omega')^2$ and $(\omega'')^2$ determine
different orientations on $M$. Choose an orthonormal basis $\{e_1,
e_2, e_3, e_4\}$ for $TM$ with $\{e_1, e_2\}$ a basis for
$T\Sigma$ such that $(\omega')^2(e_1, e_2, e_3, e_4)>0$ and
$d\mu(e_1, e_2)>0$.

Both $\omega'$ and $\omega''$ are
parallel calibrating forms and Proposition \ref{eta}
is applicable. Therefore,

\[
\begin{split}
\frac{d}{dt}\eta'=\Delta\eta'+\eta'[(h_{31k}
-h_{42k})^2
+(h_{32k}+h_{41k})^2]+(1-(\eta')^2) Ric(J_1(e_1), e_2)
\end{split}
\]

On the other hand by switching $e_3$
and $e_4$,
\[
\begin{split}
\frac{d}{dt}\eta''=\Delta\eta''+\eta''[(h_{41k}
-h_{32k})^2
+(h_{42k}+h_{31k})^2]+(1-(\eta'')^2) Ric(J_2(e_1), e_2)
\end{split}
\]
Adding these two equations and denote
$\eta'+\eta''$ by $\mu$, we get

\[
\begin{split}
\frac{d}{dt}\mu=&\Delta\mu
+\mu |A|^2
+2(\eta'-\eta'')h_{32k}h_{41k}
-2(\eta'-\eta'')h_{31k}h_{42k}\\
+&(1-(\eta')^2) Ric(J_1(e_1), e_2)
+(1-(\eta'')^2) Ric(J_2(e_1), e_2)
\end{split}
\]

Write $\mu=2 (\min\{\eta', \eta''\})+|\eta'-\eta''|$, then
$\mu\geq 2\delta+|\eta'-\eta''|$. After completing square, $\mu$
satisfies the following inequality.

\[\frac{d}{dt}\mu\geq \Delta\mu
+2\delta |A|^2-C
\]
where $-C$ is the lower bound of the
Ricci curvature of $M$, $Ric\geq -Cg$.

As before, we can isometrically embed $M$ into $\R^N$. To detect a
possible singularity at a point $(y_0, t_0)$, where $y_0\in M
\hookrightarrow \R^N$ and $t_0<\infty$, take  a ball $B$ of radius
$r$ about $y_0\in \R^N$ and $\psi$ a cut-off function as in the
proof of Proposition \ref{stat}. A similar argument yields the
following inequality:

\[
\begin{split}
&\frac{d}{dt}\int \psi (2-\mu) \rho_{y_0,t_0}\,
d\mu_t \\
\leq &C-C_\delta \int \psi \rho_{y_0,t_0} |A|^2 \,d\mu_t
\end{split}
\]
where $C_\delta$ is a constant depend on $\delta$.

Therefore $\lim_{t\rightarrow t_0} \int \psi
\rho_{y_0,t_0}(2-\mu)d\mu_t$ exists. Let $\mathcal{S}^{\lambda_i}$
be a blow-up sequence at $(y_0, t_0)$ that converges to
$\mathcal{S}^\infty$. As in the previous section we can show for a
fixed $\tau >0$,

\[\int_{-1-\tau}^{-1} \int \psi^{\lambda_j}
\rho^{\lambda_j}_{0,0} |A|^2 d\mu_s^{\lambda_j}ds\leq C(j)
\]
where $C(j)\rightarrow 0$ as $\lambda_j
\rightarrow \infty$.

Choose $\tau_j \rightarrow 0$ such that
$\frac{C(j)}{\tau_j} \rightarrow 0$ and
$s_j \in [-1-\tau_j, -1]$ so that

\[\int \psi^{\lambda_j} \rho^{\lambda_j}_{0,0} |A|^2 d\mu_{s_j}^{\lambda_j}
\leq \frac{C(j)}{\tau_j}
\]

We investigate this inequality more carefully. Notice that
$\psi^{\lambda_j}$ is supported in $B_{\lambda_j r}(0)\subset
\R^N$ and $\psi^{\lambda_j}\equiv 1$ in $B_{\frac{\lambda_j
r}{2}}(0)$. Also

\[\rho^{\lambda_j}_{0,0}(F^{\lambda_j}_{s_j})=\frac{1}{4\pi
(-s_j)}\exp(\frac{-|F^{\lambda_j}_{s_j}|^2}{4(-s_j)})\]

If we consider for any $R>0$, the ball of radius $R$,
$B_R(0)\subset \R^N$, when $j$ is large enough, we may assume
$\frac{\lambda_j r}{2}>R$ and $-1<s_j<-\frac{1}{2}$, then

\[\int \psi^{\lambda_j} \rho^{\lambda_j}_{0,0} |A|^2 d\mu_{s_j}^{\lambda_j}
\geq
\frac{1}{2\pi}\exp(\frac{-R^2}{2})\int_{\Sigma_{s_j}^{\lambda_j}\cap
B_R(0)}|A|^2 d\mu^{\lambda_j}_{s_j}\]

This implies for any compact set $K\subset \R^N$,

\[\int_{\Sigma_{s_j}^{\lambda_j}\cap K}|A|^2
d\mu_{s_j}^{\lambda_j}\rightarrow 0 \,\text{ as}\, j\rightarrow
\infty\]

Now we claim this together with the fact that $\mu$ has a positive
lower bound imply  $ \lim_{j\rightarrow \infty}\int \rho_{y_0,
t_0}d\mu_{t_0+\frac{s_j}{\lambda_j^2}} =\lim_{j\rightarrow
\infty}\int \rho_{0,0}d\mu^{\lambda_j}_{s_j}\leq 1$. We may assume
the origin in $\R^N$ is a limit point of
$\Sigma_{s_j}^{\lambda_j}$, otherwise the limit is zero.

We notice that $\omega'+\omega''$ is a parallel two form with
$\lambda_1=2$ and $\lambda_2=0$ from the last paragraph in \S 3.
Therefore the holonomy group of $M$ splits into $SO(2)\times
SO(2)$ and $M$ is locally a Riemannian product. For simplicity, we
shall assume $M$ is a product $\Sigma_1\times\Sigma_2$ such that
$\frac{1}{2}(\omega'+\omega'')$ is the volume form of $\Sigma_1$.
In fact, we can choose local coordinates $(x^1, y^1)$ on
$\Sigma_1$ and $(x^2, y^2)$ on $\Sigma_2$ so that
$\omega'=dx^1\wedge dy^1+dx^2\wedge dy^2$, $\omega''=dx^1\wedge
dy^1-dx^2\wedge dy^2$ and $\mu=2(dx^1\wedge dy^1)$.

Let $\pi_1:\Sigma_1\times \Sigma_2 \mapsto \Sigma_1$ be the
projection.  $\frac{\mu}{2}$ is in fact the Jacobian of the
projection $\pi_1$ when restricted to $\Sigma_t$ and the
restriction $\pi_1|_{\Sigma_t}$ is a covering map. Now take any
neighborhood $\Omega$ of $\pi_1(y_0)\in \Sigma_1$ and consider
$\pi_1^{-1}(\Omega) \cap \Sigma_t$. Take any component and denote
it by $\mathcal{S}_t$. $\mathcal{S}$ is an unparametrized flow.
Each $\mathcal{S}_t$ can be written as the graph of a map $u_t
:\Omega \mapsto \Sigma_2$ with uniformly bounded $|d u_t|$ since
$\mu_t$ has a uniform lower bound. Since $y_0$ is a limit point of
$\Sigma_{t_0+\frac{s_j}{\lambda_j^2}}$, by choosing $\Omega$ small
enough, we may assume the graph of $u_j=
u_{t_0+\frac{s_j}{\lambda_j^2}}$ lies in $B$.

Now we consider the parabolic blow up of the graph of $u_j$ in
$\R^N$ by $\lambda_j$. This is the graph of the map
$\widetilde{u}_j$ from $\lambda_j \Omega$ to $\lambda_j \Sigma_2$.
It corresponds to a part of $\Sigma_{s_j}^{\lambda_j}$. By the
assumption that the origin is a limit point of
$\Sigma_{s_j}^{\lambda_j}$ and $|d\widetilde{u_j}|$ is uniformly
bounded, we may assume $\widetilde{u}_j \rightarrow
\widetilde{u}_\infty$ in $C^\alpha$ on compact sets.
$\widetilde{u}_\infty$ is an entire graph defined on $\R^2$.

Other the other hand,

\[|A|_j\leq |\nabla d\widetilde{u}_j|\leq
(\sqrt{1+|d\widetilde{u}_j|^2})^3|A|_j\] where $|A|_j$ is the norm
of the second fundamental form of $\mathcal{S}_{s_j}^{\lambda_j}$
and $|\nabla d\widetilde{u}_j|$ is the norm of the covariant
derivatives of $d\widetilde{u}_j$. Now we identify $\Omega$ with
an open set in $\R^2$. Therefore for any $B_\rho \subset \R^2$,
 $\widetilde{u}_j$ satisfies

\[ |D \widetilde{u}_j|\leq C, \int_{B_\rho}
|D^2  \widetilde{u}_j|^2 \rightarrow 0\] where $D \widetilde{u}_j$
and $D^2  \widetilde{u}_j$ are the usual derivatives with respect
to coordinate variables on $\R^2$. Denote $v_j=\frac{\partial
\widetilde{u}_j} {\partial x_k}$, then $|v_j|\leq C$ and
$\int_{B_\rho}|Dv_j|^2 \rightarrow 0$.

Let $c_j=\frac{1}{Vol(B_\rho)}\int v_j$, then we can choose a
convergent subsequence $c_j\rightarrow c$. By Poincare inequality,

\[\int|v_j-c_j|^2
\leq \lambda \int |Dv_j|^2 \rightarrow 0\]

Therefore $\frac{\partial  \widetilde{u}_j}{\partial x^k}
\rightarrow c_k$ in $L^2$. Since we may
assume $ \widetilde{u}_j \rightarrow
 \widetilde{u}_\infty$ in
$C^\alpha\cap W^{1,2}_{loc}$, this implies
$\mathcal{S}^{\lambda_j}_{s_j} \rightarrow
\mathcal{S}^\infty_{-1}$ as Radon measures and
$\mathcal{S}^\infty_{-1}$ is the graph of a linear function.
Therefore

\[\lim_{j\rightarrow \infty}\int \rho_{0,0}
d\mu_{s_j}^{\lambda_j} =\int \rho_{0,0} d\mu^\infty_{-1}=1\] By
White's theorem again, we have regularity at the point $y_0,
t_0)$.

\end{proof}

Now we prove Theorem B.

\vskip 5pt
\noindent {\bf Theorem B}
\hskip 3pt{\it Let $M$ be an oriented four-dimensional
Einstein manifold with two parallel calibrating
forms $\omega', \omega''$ such that $\omega'$
is self-dual and $\omega''$ is anti-self-dual.
If $\Sigma$ is a compact oriented surface
immersed in $M$ such that $*\omega',
*\omega''>0$ on $\Sigma$. Then the mean
curvature flow of $\Sigma $ exists
smoothly for all time.}
\vskip 8pt

\begin{proof}
$*\omega'$ and
$*\omega''$ have positive lower bound for
any finite time by equation (\ref{eta2}),
therefore the assumption in Proposition
\ref{long} is satisfied.
\end{proof}

\section{Convergence at infinity}

In this section we study the convergence of the
mean curvature flow at infinity. The key point
is to show uniform boundedness of $|A|^2$
in space and time.
We first compute the evolution of
the second fundamental form.
Let $\Sigma \rightarrow M^n$ be an isometric
immersion.
We choose an orthonormal basis $\{e_i\}$ for
$T\Sigma$ and
$\{e_\alpha\}$ for $N\Sigma$.
Recall the convention for indexes are
$A, B, C \cdots =1\cdots n$, $i,j,k \cdots$ for tangent
indexes, and
$\alpha, \beta, \gamma \cdots$ for normal indexes. Now denote
the coefficient of the second fundamental
form by $h_{\alpha ij}=<A(\pai,\pj), e_\alpha>$.
The covariant derivative of $A$ is
defined as

\[(\barbla_{\pl}A)(\pai, \pj)
=(\barbla_{\pl} A(\pai, \pj))^N
-A((\barbla_{\pl} \pai)^T, \pj)
-A(\pai, (\barbla_{\pl}\pj)^T)
\]

We denote

\[h_{\alpha ij,k}=<(\barbla_{\pk}A)(\pai, \pj)
,e_\alpha>\]

and
\[h_{\alpha ij,kl}=<(\barbla_{\pl}\barbla
_{\pk}A)(\pai, \pj)
,e_\alpha>\]

Let $\Delta h_{\alpha ij}
=g^{kl}h_{\alpha ij,kl}$ be the Laplacian
of $h_{\alpha ij}$.

\begin{pro}
For a mean curvature flow $F:\Sigma\times
[0,t_0)\rightarrow M$ of any dimension,
the second fundamental form $h_{\alpha ij}$
satisfies the following equation.

\begin{equation}\label{h}
\begin{split}
\frac{d}{dt}h_{\alpha ij}
&=\Delta h_{\alpha ij}
+(\barbla_{\pk}K)_{\alpha ijk}
+(\barbla_{\pj}K)_{\alpha kik}\\
&-2K_{lijk}h_{\alpha lk}
+2K_{\alpha \beta jk}h_{\beta ik}
+2K_{\alpha \beta ik}h_{\beta jk}\\
&-K_{lkik}h_{\alpha lj}
-K_{lkjk}h_{\alpha li}
+K_{\alpha k\beta k}h_{\beta ij}\\
&-h_{\alpha im}(h_{\gamma mj}h_{\gamma}
-h_{\gamma mk}h_{\gamma jk})\\
&-h_{\alpha mk}(h_{\gamma mj}h_{\gamma ik}
-h_{\gamma mk}h_{\gamma ij})\\
&-h_{\beta ik}(h_{\beta lj}h_{\alpha lk}
-h_{\beta lk}h_{\alpha lj})\\
&-h_{\alpha jk}h_{\beta ik}h_\beta
+h_{\beta ij}<e_\beta, \barbla_{H} e_\alpha>
\end{split}
\end{equation}
where $K_{ABCD}$ is the curvature tensor
and $\barbla$ is the covariant derivative of
$M$.

In particular, $|A|^2$ satisfies the following
equation along the mean curvature flow.
\begin{equation}\label{|A|^2}
\begin{split}
\frac{d}{dt}|A|^2
&=\Delta |A|^2 -2|\nabla A|^2
+2[(\barbla_{\pk}K)_{\alpha ijk}
+(\barbla_{\pj}K)_{\alpha kik}]h_{\alpha ij}\\
&-4K_{lijk}h_{\alpha lk}h_{\alpha ij}
+8K_{\alpha \beta jk}h_{\beta ik}h_{\alpha ij}
-4K_{lkik}h_{\alpha lj}h_{\alpha ij}
+2K_{\alpha k\beta k}h_{\beta ij}h_{\alpha ij}\\
&+2\sum_{\alpha,\gamma, i,m}
(\sum_k h_{\alpha ik}h_{\gamma mk}
-h_{\alpha mk}h_{\gamma ik})^2
+2\sum_{i,j,m,k}(\sum_{\alpha} h_{\alpha ij}
h_{\alpha mk})^2
\end{split}
\end{equation}
\end{pro}

\begin{proof}

We first derive equation (\ref{|A|^2})
from equation (\ref{h}).

Since $|A|^2=g^{ik}g^{jl}h_{\alpha ij}
h_{\alpha kl}$, calculate using a normal coordinate
system near a point $p$ we have

\[\frac{d}{dt}|A|^2
=2(\frac{d}{dt}g^{ik})h_{\alpha ij}
h_{\alpha kj}
+2(\frac{d}{dt} h_{\alpha ij})h_{\alpha ij}
\]

Recall $\frac{d}{dt} g_{ik}
=2h_\beta h_{\beta ik}$ and plug in equation (
\ref{h}) to get

\[
\begin{split}
\frac{d}{dt}|A|^2
&=4h_\beta h_{\beta ik}h_{\alpha ij}h_{\alpha kj}\\
&+2h_{\alpha ij}[\Delta h_{\alpha ij}
+(\barbla_{\pk}K)_{\alpha ijk}
+(\barbla_{\pj}K)_{\alpha kik}\\
&-2K_{lijk}h_{\alpha lk}
+2K_{\alpha \beta jk}h_{\beta ik}
+2K_{\alpha \beta ik}h_{\beta jk}\\
&-K_{lkik}h_{\alpha lj}
-K_{lkjk}h_{\alpha li}
+K_{\alpha k\beta k}h_{\beta ij}\\
&-h_{\alpha im}(h_{\gamma mj}h_{\gamma}
-h_{\gamma mk}h_{\gamma jk})\\
&-h_{\alpha mk}(h_{\gamma mj}h_{\gamma ik}
-h_{\gamma mk}h_{\gamma ij})\\
&-h_{\beta ik}(h_{\beta lj}h_{\alpha lk}
-h_{\beta lk}h_{\alpha lj})\\
&-h_{\alpha jk}h_{\beta ik}h_\beta
+h_{\beta ij}<e_\beta, \barbla_{H} e_\alpha>]
\end{split}
\]
The first term on the right hand side
$4h_\beta h_{\beta ik}h_{\alpha ij}h_{\alpha kj}$
cancels with two later terms. They are
so-called "metric" terms and vanish if we
choose a orthonormal frame in our computation.

The last term on the right
 hand side
 $2 h_{\alpha ij}h_{\beta ij}
 <e_\beta, \barbla_{H} e_\alpha>$ is zero by
 symmetry.

Now use
\[\Delta h_{\alpha ij}^2
=2|\nabla A|^2 +2h_{\alpha ij}\Delta h_{\alpha
ij}\]
Therefore we get

\[
\begin{split}
\frac{d}{dt}|A|^2
&=\Delta |A|^2-2|\nabla A|^2\\
&+2h_{\alpha ij}[(\barbla_{\pk}K)_{\alpha ijk}
+(\barbla_{\pj}K)_{\alpha kik}\\
&-2K_{lijk}h_{\alpha lk}
+2K_{\alpha \beta jk}h_{\beta ik}
+2K_{\alpha \beta ik}h_{\beta jk}\\
&-K_{lkik}h_{\alpha lj}
-K_{lkjk}h_{\alpha li}
+K_{\alpha k\beta k}h_{\beta ij}\\
&+h_{\alpha im}h_{\gamma mk}h_{\gamma jk}-h_{\alpha mk}(h_{\gamma mj}h_{\gamma ik}
-h_{\gamma mk}h_{\gamma ij})-h_{\beta ik}(h_{\beta lj}h_{\alpha lk}
-h_{\beta lk}h_{\alpha lj})]
\end{split}
\]

The fourth order terms can be calculated
as the following

\[
\begin{split}
&h_{\alpha ij}h_{\alpha im}
h_{\gamma mk}h_{\gamma jk}
-h_{\alpha ij}h_{\alpha mk}(h_{\gamma mj}h_{\gamma ik}
-h_{\gamma mk}h_{\gamma ij})
-h_{\alpha ij}h_{\beta ik}(h_{\beta lj}h_{\alpha lk}
-h_{\beta lk}h_{\alpha lj})\\
&=2h_{\alpha ij}h_{\alpha im}h_{\gamma mk}h_{\gamma kj}
-2h_{\alpha ij}h_{\alpha mk}h_{\gamma mj}h_{\gamma ik}
+h_{\alpha ij}h_{\alpha mk}h_{\gamma mk}h_{\gamma ij}\\
\end{split}
\]

The first two terms can be completed to square.

\begin{equation}
\begin{split}
&2h_{\alpha ij}h_{\alpha im}
h_{\gamma mk}h_{\gamma kj}
-2h_{\alpha ij}h_{\alpha mk}
h_{\gamma mj}h_{\gamma ik}\\
&=2h_{\alpha ij}h_{\alpha ik}
h_{\gamma mk}h_{\gamma mj}
-2h_{\alpha ij}h_{\alpha mk}
h_{\gamma mj}h_{\gamma ik}\\
&=2h_{\alpha ij}h_{\gamma mj}
(h_{\alpha ik}h_{\gamma mk}
-h_{\alpha mk}h_{\gamma ik})\\
&=h_{\alpha ij}h_{\gamma mj}
(h_{\alpha ik}h_{\gamma mk}
-h_{\alpha mk}h_{\gamma ik})
+h_{\alpha mj}h_{\gamma ij}
(h_{\alpha mk}h_{\gamma ik}
-h_{\alpha ik}h_{\gamma mk})\\
&=\sum_{\alpha,\gamma, i,m}
(\sum_k h_{\alpha ik}h_{\gamma mk}
-h_{\alpha mk}h_{\gamma ik})^2
\end{split}
\end{equation}

Now we calculate the equation (\ref{h}).
First the Laplacian of $h_{\alpha ij}$
is the following.

\begin{equation}\label{lap}
\begin{split}
\Delta h_{\alpha ij}
=&h_{\alpha ,ij}-(\nabla_{\pk}K)_{\alpha ijk}
-(\nabla_{\pj}K)_{\alpha kik}\\
&+2K_{lijk}h_{\alpha lk}
-2K_{\alpha \beta jk}h_{\beta ik}
-2K_{\alpha \beta ik}h_{\beta jk}
-K_{\alpha ij\beta}h_{\beta}\\
&+K_{lkik}h_{\alpha lj}
+K_{lkjk}h_{\alpha li}
-K_{\alpha k\beta k}h_{\beta ij}\\
&+h_{\alpha im}(h_{\gamma mj}h_{\gamma}
-h_{\gamma mk}h_{\gamma jk})\\
&+h_{\alpha mk}(h_{\gamma mj}h_{\gamma ik}
-h_{\gamma mk}h_{\gamma ij})\\
&+h_{\beta ik}(h_{\beta lj}h_{\alpha lk}
-h_{\beta lk}h_{\alpha lj})\\
\end{split}
\end{equation}
where $h_{\alpha ,ij}
=<\barbla_{\pj}^N \barbla_{\pai}^N H, e_\alpha>$.

In codimension one case, this equation
reduces to
\[
\begin{split}
\Delta h_{ij}
=&H_{,ij}-(\nabla_{\pk}K)_{N ijk}
-(\nabla_{pj}K)_{N kik}\\
&+2K_{lijk}h_{lk}
-K_{N ij n+1}H\\
&+K_{lkik}h_{ lj}
+K_{lkjk}h_{li}
-K_{N k N k}h_{ ij}\\
&+h_{im}h_{ mj}H
-h_{mk}^2 h_{ij}\\
\end{split}
\]
This recovers equation (1.20) in \cite{ssy}.

The equation (\ref{lap}) is computed using the Codazzi equation and the commutation
formula.
\[
\begin{split}
&h_{\alpha kj,i}=h_{\alpha ki,j}+K_{\alpha kij}\\
&h_{\alpha ij,kl}=h_{\alpha ij,lk}
-h_{\alpha im}R_{mjlk}
-h_{\alpha mj}R_{milk}
-h_{\beta ij}R_{\beta\alpha lk}
\end{split}
\]

where $R_{mjlk}$ is the curvature of $T\Sigma$
and $R_{\alpha \beta lk}$ is the curvature
of $N\Sigma$.

We start the computation with $h_{\alpha, ij}
=h_{\alpha kk,ij}$.

\[
\begin{split}
&h_{\alpha kk,ij}\\
&=(h_{\alpha ki,k}+K_{\alpha kik})_j\\
&=h_{\alpha ki,kj}+K_{\alpha kik,j}\\
&=h_{\alpha ik,jk}
-h_{\alpha im}R_{mkjk}
-h_{\alpha mk}R_{mijk}
-h_{\beta ik}R_{\beta\alpha jk}+K_{\alpha kik,j}\\
&=(h_{\alpha ij,k}+K_{\alpha ijk})_k
-h_{\alpha im}R_{mkjk}
-h_{\alpha mk}R_{mijk}
-h_{\beta ik}R_{\beta\alpha jk}+K_{\alpha kik,j}\\
\end{split}
\]

By the Gauss and Ricci equation, we have
\[
\begin{split}
&R_{mkjk}=K_{mkjk}+h_{\gamma mj}h_{\gamma kk}
-h_{\gamma mk}h_{\gamma kj}\\
&R_{mijk}=K_{mijk}+h_{\gamma mj}h_{\gamma ik}
-h_{\gamma mk}h_{\gamma ij}\\
&R_{\beta\alpha jk}=K_{\beta \alpha jk}
+h_{\beta lj}h_{\alpha lk}
-h_{\beta lk}h_{\alpha lj}\\
\end{split}
\]

Therefore

\[
\begin{split}
\Delta h_{\alpha ij}
&=h_{\alpha ,ij}-K_{\alpha ijk,k}-K_{\alpha kik,j}
+h_{\alpha im}K_{mkjk}
+h_{\alpha mk}K_{mijk}
+h_{\beta ik}K_{\beta\alpha jk}\\
&+h_{\alpha im}(h_{\gamma mj}h_{\gamma}
-h_{\gamma mk}h_{\gamma jk})\\
&+h_{\alpha mk}(h_{\gamma mj}h_{\gamma ik}
-h_{\gamma mk}h_{\gamma ij})\\
&+h_{\beta ik}(h_{\beta lj}h_{\alpha lk}
-h_{\beta lk}h_{\alpha lj})\\
\end{split}
\]

The covariant derivative term can
be calculated
as the following.
\[
\begin{split}
&K_{\alpha ijk,k}\\
&=(\barbla_{\pk}K)_{\alpha ijk}-K_{lijk}h_{\alpha lk}
+K_{\alpha \beta jk}h_{\beta ik}
+K_{\alpha i\beta k}h_{\beta jk}
+K_{\alpha ij\beta}h_{\beta kk}\\
&K_{\alpha kik,j}\\
&=(\barbla_{\pj}K)_{\alpha kik}-K_{lkik}h_{\alpha lj}
+K_{\alpha\beta ik}h_{\beta kj}
+K_{\alpha k\beta k}h_{\beta ij}
+K_{\alpha ki\beta }h_{\beta kj}
\end{split}
\]

Note that $K_{\alpha ijk}$ is considered as a section of the
bundle $N\Sigma \otimes T\Sigma\otimes T\Sigma \otimes T\Sigma$ in
taking covariant derivatives . We collect all the embient
curvature term and use the first Bianchi identity $K_{\alpha
i\beta k} +K_{\alpha k i\beta} =-K_{\alpha \beta ki}$ to get
equation (\ref{lap}).

Next we calculate the equation for
$h_{\alpha ij}=<\barbla_{\pj}\pai, e_\alpha>$.
\[
\begin{split}
\frac{d}{dt}h_{\alpha ij}
&=<\barbla_H \barbla_{\pj}\pai, e_\alpha>
+<\barbla_{\pj}\pai, \barbla_H {e_\alpha}>\\
&=<\barbla_{\pj}\barbla_H \pai, e_\alpha>
-<K(H,\pj)\pai, e_\alpha>
+<\barbla_{\pj}\pai, \barbla_{H} e_\alpha>
\end{split}
\]

By breaking
 $\barbla_{\pj}\barbla_{\pai} H$ into normal and tangent parts,
we get
\[
\begin{split}
<\barbla_{\pj}\barbla_{\pai} H, e_\alpha>
&=<\barbla_{\pj}[(\barbla_{\pai} H)^T
+(\barbla_{\pai} H)^N], e_\alpha>\\
&=<\barbla_{\pj}^N \barbla_{\pai}^N H, e_\alpha>
-<(\barbla_{\pai}H)^T, \barbla_{\pj} e_\alpha>
\end{split}
\]

Therefore,

\[
\begin{split}
\frac{d}{dt}h_{\alpha ij}
=h_{\alpha, ij}-h_\beta K_{\beta ji\alpha}
-<(\barbla_{\pai}H)^T, \barbla_{\pj} e_\alpha>
+<\barbla_{\pj} \pai, \barbla_H e_\alpha>
\end{split}
\]
where
$h_{\alpha, ij}=<\barbla_{\pj}^N
\barbla_{\pai}^N H, e_\alpha>$.

The term $<(\barbla_{\pai}H)^T,
\barbla_{\pj} e_\alpha>$ is equal
to $h_\beta h_{\beta ik}h_{\alpha jk}$.
Also since we choose a normal coordinate
in our computation, $(\barbla_{\pai}\pj)^T
=0$ and
$<\barbla_{\pai} \pj, \barbla_H e_\alpha>
=h_{\beta ij}<e_\beta, \barbla_H e_\alpha>$.

\begin{equation}\label{ddt}
\begin{split}
\frac{d}{dt}h_{\alpha ij}
=h_{\alpha, ij}-h_{\beta} h_{\beta ik}
h_{\alpha jk} -h_\beta K_{\beta ji\alpha}
+h_{\beta ij} <e_\beta, \barbla_H e_\alpha>
\end{split}
\end{equation}

Combine equation (\ref{lap}) and
(\ref{ddt}), we get the parabolic equation
for $h_{\alpha ij}$.

\end{proof}

The following proposition provides a uniform bound
of the second fundamental form when
$*\omega'$ and $*\omega''$ are both close
to one.

\begin{pro}\label{conv}
Let $M$ be a compact four-dimensional
manifold with
bounded geometry.  Let $\omega'$ and $\omega''$
be two parallel calibrating forms such that
$\omega'\wedge \omega'$ and $\omega''
\wedge \omega'' $ determine opposite
orientation for $M$. Let $\Sigma$ be
an oriented immersed surface in $M$.
There exists a constant $1>\epsilon>0$  such
that if $*\omega'>1-\epsilon$ and
$*\omega''>1-\epsilon$
on $F_t(\Sigma)$ for  $t\in [0,T]$, then the
norm of the second fundamental form of
$F_t(\Sigma)$ is uniformly bounded in $[0,T]$.
\end{pro}

\begin{proof}
The fourth order term in equation (\ref{|A|^2})
can be calculated explicitly in the
four-dimensional case.

\begin{equation}
\begin{split}
&\sum_{\alpha, \gamma, i,m}(\sum_k h_{\alpha ik}h_{\gamma mk}
-h_{\alpha mk}h_{\gamma ik})^2\\
&=\sum_{i,m}(\sum_k h_{3ik}h_{4mk}-h_{3mk}h_{4ik})^2
+\sum_{i,m}(\sum_k h_{4ik}h_{3mk}-h_{4mk}h_{3ik})^2\\
&=2\sum_{i,m}(\sum_k h_{3ik}h_{4mk}-h_{3mk}h_{4ik})^2\\
&=2[(\sum_kh_{31k}h_{42k}-h_{32k}h_{41k})^2
+(\sum_k h_{41k}h_{32k}-h_{42k}h_{31k})^2]\\
&=4(\sum_k h_{31k}h_{42k}-h_{32k}h_{41k})^2\\
&=4[\frac{1}{2}(h_{31k}+h_{42k})^2
+\frac{1}{2}(h_{32k}-h_{41k})^2
-\frac{1}{2}|A|^2]^2\\
&=[(h_{31k}+h_{42k})^2+(h_{32k}-h_{41k})^2
-|A|^2]^2
\end{split}
\end{equation}

By Schwarz inequality,
\[
\sum_{i,j, m,k}(\sum_\alpha
h_{\alpha ij} h_{\alpha mk})^2
\leq \sum_{i,j, m, k}(\sum_\alpha h_{\alpha ij}^2)
(\sum_\alpha h_{\alpha mk}^2)
\leq |A|^4
\]

Therefore
\begin{equation}\label{A2}
\begin{split}
\frac{d}{dt} |A|^2
&\leq \Delta |A|^2 -2|\nabla A|^2+4|A|^4+K_1 |A|^2+K_2\\
\end{split}
\end{equation}
where $K_1$ and $K_2$ are constants that
depend on the curvature tensor
and covariant derivatives of the curvature
tensor of $M$.

Again we consider $\mu=\eta'+\eta''$.
Since $\eta'\geq 1-\epsilon$ and $\eta''\geq
1-\epsilon$, we have $\mu\geq 2-2\epsilon$
and $|\eta'-\eta''|\leq \epsilon
\leq \frac{\epsilon}{2-2\epsilon}\mu$.

\[
\begin{split}
\frac{d}{dt}\mu
&\geq \Delta\mu
+\mu |A|^2
+2(\eta'-\eta'')h_{32k}h_{41k}
-2(\eta'-\eta'')h_{31k}h_{42k}-C\mu\\
& \geq \Delta\mu
+\mu (|A|^2
-2\frac{\epsilon}{2-2\epsilon}|h_{32k}h_{41k}|
-2\frac{\epsilon}{2-2\epsilon}|h_{31k}h_{42k}|)-C\mu\\
&\geq \Delta\mu
+\mu [\frac{2-3\epsilon}{2-2\epsilon}|A|^2
+\frac{\epsilon}{2-2\epsilon}
(|A|^2-2|h_{32k}h_{41k}|-2|h_{31k}h_{42k}|)]-C\mu\\
\end{split}
\]
The term $(|A|^ 2-2|h_{32k}h_{41k}|-2|h_{31k}h_{42k}|)$
is a complete square and thus nonnegative.

Therefore
\[
\frac{d}{dt} \mu \geq \Delta\mu +c_1(\epsilon)\mu |A|^2-C\mu
\]
where $c_1(\epsilon)=\frac{2-3\epsilon}
{2-2\epsilon}$, notice that $c_1(\epsilon)
\rightarrow 1$ as $\epsilon\rightarrow 0$.

Let $p>1$ be an integer to be determined, we
calculate the equation for $\mu^p$.

\[
\frac{d}{dt}\mu^p=p\mu^{p-1}\frac{d}{dt}\mu
\geq p\mu^{p-1}(\Delta \mu+c_1(\epsilon)\mu
|A|^2-C\mu)
\]
Use the identity
$\Delta \mu^p=p(p-1)\mu^{p-2}
|\nabla \mu|^2+p\mu^{p-1}\Delta\mu$, the
differential inequality for $\mu^p$ becomes

\[\frac{d}{dt}\mu^p\geq\Delta \mu^p-p(p-1)
\mu^{p-2}|\nabla \mu|^2
+p\,c_1(\epsilon)\mu^p|A|^2-Cp\mu^p
\]

Now we estimate the term $|\nabla \mu|^2$

We calculate as in equation (\ref{deta})

\[
\begin{aligned}
|\nabla\eta'|^2&\leq
2({1-(\eta')^2})(h_{41k}^2+h_{32k}^2)\\
|\nabla \eta''|^2&\leq 2({1-(\eta'')^2})
(h_{31k}^2+h_{42k}^2)
\end{aligned}\]

Since $\eta'\geq 1-\epsilon$, $1-(\eta')^2
\leq 1-(1-\epsilon)^2\leq 2\epsilon$.
Therefore
\[|\nabla \mu|^2
\leq 2(|\nabla\eta'|^2+|\nabla\eta''|^2)
\leq 8\epsilon |A|^2
\leq\frac{8\epsilon}{(2-2\epsilon)^2}
\mu^2|A|^2
\]
Denote $\frac{8\epsilon}{(2-2\epsilon)^2}=c_2(
\epsilon)$.

Thus
\[\frac{d}{dt}\mu^p\geq\Delta \mu^p+p[
c_1(\epsilon)-(p-1)c_2(\epsilon)]
\mu^p|A|^2-Cp\mu^p
\]

Plug in $f=|A|^2$ and $g=\mu^p $ in
the identity

\[\frac{d}{dt}(\frac{f}{g})
=\Delta(\frac{f}{g})
+2\frac{\nabla g}{g}\cdot \nabla(\frac{f}{g})
+\frac{1}{g^2}
[(\frac{d}{dt}f-\Delta f)g
-(\frac{d}{dt}g-\Delta g)f]
\]

Therefore
\[
\begin{split}
\frac{d}{dt}(\frac{|A|^2}{\mu^p})
\leq &\,\, \Delta (\frac{|A|^2}{\mu^p})
+2\nabla (\frac{|A|^2}{\mu^p})\cdot
\frac{\nabla \mu^p}{\mu^p}\\
+&\frac{1}{\mu^{2p}}
\{[-2|\nabla A|^2
+4|A|^4+K_1|A|^2+K_2]\mu^p\\
&\,\,\,\,\,-[p(c_1(\epsilon)-(p-1)c_2(\epsilon))\mu^p
|A|^2-Cp\mu^p]|A|^2\}
\end{split}
\]

The last term is less than

\[[4-p(c_1(\epsilon)-(p-1)c_2(\epsilon))]
\frac{|A|^4}{\mu^p}+(K_1+Cp)
\frac{|A|^2}{\mu^p}+K_2 \frac{1}{\mu^p}
\]

Recall that $c_1(\epsilon)\rightarrow 1$ and
$c_2(\epsilon)\rightarrow 0$ as
$\epsilon \rightarrow 0$.

Choose $p$ large enough and then $\epsilon$
small enough so that
$(2-2\epsilon)^p
[4-p(c_1(\epsilon)-(p-1)c_2(\epsilon))]\leq
-C_1
$ for some $C_1>0$

Then
\[
[4-p(c_1(\epsilon)-(p-1)c_2(\epsilon))]
\frac{|A|^4}{\mu^p}
=[4-p(c_1(\epsilon)-(p-1)c_2(\epsilon))]\mu^p
\frac{|A|^4}{\mu^{2p}}
\leq- C_1 \frac{|A|^4}{\mu^{2p}}
\]

Denote $f=\frac{|A|^2}{\mu^{p}}$, then $f$ satisfies

\[
\begin{split}
\frac{d}{dt}f
\leq \, \Delta f
+ V\cdot \nabla f
-C_1 f^2 +C_2 f+C_3
\end{split}
\]

Now we apply the maximum principle for parabolic equations and
conclude the $\frac{|A|^2}{\mu}$ is uniformly bounded, thus $|A|^2
$ is also bounded.

\end{proof}

Now we prove Theorem C.

\vskip 5pt
\noindent {\bf Theorem C}
\hskip 3pt{\it Under the same assumption as
in Theorem B. When $M$ has non-negative
curvature, there exists a constant
$\epsilon>0$ such that if $\Sigma$ is a compact oriented surface
immersed in $M$ with $*\omega',
*\omega''>1-\epsilon $ on $\Sigma$, the mean
curvature flow of $\Sigma $ converges
smoothly to a totally geodesic
surface at infinity.}
\vskip 8pt
\begin{proof} In these cases, $*\omega'$ and $*\omega''$
are both non-decreasing. Proposition
\ref{conv} is applicable and $|A|^2$
is uniformly bounded in space and time.

Integrating equation (\ref{A2}) and we
see
\begin{equation}\label{bb}
\frac{d}{dt}\int_{\Sigma_t}
|A|^2 d\mu_t\leq C
\end{equation}
 Recall in this case
$\frac{d}{dt}\mu  \geq \Delta\mu
+c_1(\epsilon) \mu |A|^2$ and $\eta$ has
a positive lower bound, thus

\begin{equation}\label{cc}
\int_0^\infty \int_{\Sigma_t}
|A|^2 d\mu_t dt\leq \infty
\end{equation}

Equation (\ref{bb}) and (\ref{cc})
together
implies

\[\int_{\Sigma_t}|A|^2d\mu_t \rightarrow
0\]

By the small $\epsilon$ regularity theorem in
\cite{il1}, $\sup_{\Sigma_t}|A|^2\rightarrow 0$ uniformly
as $t\rightarrow \infty$.

 Since the mean curvature
flow is a gradient flow and the metrics are
analytic, by the theorem of Simon \cite{si},
we get convergence at infinity.

\end{proof}

\section{Applications}

In the following we apply the previous theorem to
the case when $M$ is a product $S^2\times
S^2 $. Denote their K\"ahler forms by
$\omega_1$ and $\omega_2$ respectively.

Let $\omega'=\omega_1+\omega_2$ and $\omega''=\omega_1-\omega_2$,
then $(\omega')^2=2\omega_1\wedge \omega_2$ and
$(\omega'')^2=-2\omega_1\wedge \omega_2$ determine opposite
orientations on $M$. Both $\omega'$ and $\omega''$ are parallel
calibrating form and they define integrable almost complex
structures with opposite orientations.

\vskip 10pt \noindent {\bf Theorem D} \hskip 3pt{\it Let $M=(S^2,
\omega_1)\times (S^2,\omega_2)$. If $\Sigma$ is a compact oriented
surface embedded in $M$ such that $*\omega_1\geq|*\omega_2|$ and
the strict inequality holds at at least one point. Then the mean
curvature flow of $\Sigma$ exists for all time and converges
smoothly to a $S^2\times \{p\}$.} \vskip 10pt

\begin{proof}
We notice the statement is a little bit different from the one
given in the introduction. The difference is resolved by the
considering the maximum principle for $\eta'=*\omega'$ and
$\eta''=*\omega''$. It is not hard to see the assumption implies
$*\omega_1>|*\omega_2|$ holds everywhere at a later time. By the
equation of $\eta' $ and $\eta''$,

\[
\begin{split}
\frac{d}{dt}\eta'=\Delta\eta'+\eta'[(h_{31k}
-h_{42k})^2
+(h_{32k}+h_{41k})^2]+\eta'(1-(\eta')^2)
\end{split}
\]
\[
\begin{split}
\frac{d}{dt}\eta''=\Delta\eta''+\eta''[(h_{41k}
-h_{32k})^2
+(h_{42k}+h_{31k})^2]+\eta''(1-(\eta'')^2)
\end{split}
\]
 we
see that as $t\rightarrow \infty$, they both approach $1$. By
Theorem A, we have existence for all time.  Also the assumption on
Proposition \ref{conv} is satisfied and the second fundamental
form is uniformly bounded in space and time.  Since the mean
curvature flow is a gradient flow and the metrics are analytic, we
can apply Simon's theorem \cite{si}  to conclude convergence at
infinity. The limiting submanifold has $*\omega_1=1$ identically
and thus is of the form $S^2\times p$.

\end{proof}

Corollary D now follows from this since the
condition $*\omega_1>|*\omega_2|$ on the
graph of a map $f$ is equivalent to the
Jacobian of $f$ being less than one.

We conclude this section by the following two remarks,
\begin{rem}
When $M$ is locally a product of two Riemann surfaces of
nonpositive curvature, the method in $\cite{mu1}$ can be used to
prove uniform convergence of the flow.  The limit is totally
geodesic and the corresponding map converges to one ranged in a
lower dimensional submanifold. Notice that the convergence in
Theorem C is a stronger smooth convergence under the closeness
assumption. Such results are generalized to arbitrary dimension
and codimension in \cite{mu}.
\end{rem}

\begin{rem}
The case when $*\omega_1>0$ and $*\omega_2=0$ corresponds to
$\Sigma$ is Lagrangian surface with respect to the symplectic form
$\omega_2$. If $\Sigma $ is the graph of a map $f$, then $f$ is
indeed an area preserving diffeomorphism. This  case and the
application to the structure of the diffeomorphism groups of
compact Riemann surfaces are discussed in \cite{mu1}.

\end{rem}

\end{document}